\documentclass[12pt,a4paper]{article}

\usepackage{amssymb}
\usepackage{amsthm}
\usepackage{dsfont}
\usepackage{amsmath}
\usepackage{graphicx}
\usepackage[T1]{fontenc}
\usepackage[utf8]{inputenc}
\usepackage{tikz}
\usepackage{pxfonts}

\swapnumbers
\newtheorem{proposition}{Proposition}
\newtheorem{example}[proposition]{Example}
\newtheorem{lemma}[proposition]{Lemma}

\newtheorem{theorem}[proposition]{Theorem}
\newtheorem{remark}[proposition]{Remark}

\newtheorem{observation}[proposition]{Observation}

\numberwithin{proposition}{section}

\def\prf{\noindent{\it Proof. }}
\def\qed{\hfill$\square$}

\def\mC{\mathcal{C}}
\def\mF{\mathcal{F}}
\def\mP{\mathcal{P}}
\def\mM{\mathcal{M}}
\def\mN{\mathcal{N}}

\def\supp{\operatorname{supp}}

\title{Non-exchangeability of copulas arising from shock models}
\author{Damjana Kokol Bukov\v sek, Toma\v z Ko\v sir, Bla\v z Moj\v skerc, and Matja\v z Omladi\v c\footnote{All four authors acknowledge financial support of the Slovenian Research
Agency (research core funding No. P1-0222).}}%
\begin{document}

\maketitle

%\begin{center}

{\bf %\large
\noindent
Abstract}. Copula is a useful tool that captures the dependence structure among random variables. In practice, it is an important question which copula to choose depending on the given data and stochastic assumptions on the model in order to achieve an appropriate interpretation of the data at hand. This paper intends to help a practitioner to make a better decision about that. We concentrate on the study of the lack of exchangeability, a copulas' attribute closely studied only recently. The main non-exchangeability measure $\mu_\infty$ for a family of copulas is the supremum of the differences $|C(x,y)-C(y,x)|$ over all $(x,y)$ and all copulas $C$ in the family. We give the sharp bound of $\mu_{\infty}$ for the families of Marshall copulas, maxmin and reflected maxmin copulas (i.e.\ the main shock-model based copulas) as well as the families of positively and of negatively quadrant dependent copulas. A major contribution of this paper is also exact calculation of the maximal asymmetry function on each of the particular families of copulas. When restricted to special families of copulas considered, it helps us finding the sharp bound of $\mu_\infty$ for each of the given families. And even more importantly, it helps us giving a stochastic interpretation of the extremal copulas and examples of shock models where the maximal asymmetry is attained.

%\end{center}

\medskip

{\small \noindent


  }\\[,3 cm]

\medskip

\section{ Introduction }

Copulas are mathematical objects that capture the dependence structure among random variables. Since they were introduced by {A.\ Sklar in 1959} they have gained a lot of popularity and applications in several fields, e.g., in finance, insurance and reliability theory. Through them we study measures of dependence and build families of distributions with given margins.

An important class of copulas for applications are those arising from shock models: Marshall copulas, maxmin copulas, and reflected maxmin copulas (RMM for short). These copulas have a long history starting with \cite{MaOl} and \cite{M} and going up to \cite{KKMO}, say, where an extensive overview of these models is given together with an appropriate bibliography. A comprehensive list of references of concrete applications of shock-based copulas would be too long to present here, so let us limit ourselves to four of them, relatively recent ones and in quite different fields: \cite{LiMcNe,ChMu,DuOk,Hu}. Note that our investigations are not only of a theoretical interest, but also of a practical impact in the construction of statistical models (in a parametric as well as non-parametric context). When choosing the right copula for the data at hand a key point is to determine the family that describes the phenomenon behind the data at the best. In this respect, a better choice of the copulas could be obtained when the information about the non-symmetry of the data (measured, for instance, in a nonparametric way by means of the empirical copula) is also taken into account (cf.\ \cite[Section 5]{DKMS} where this point is further discussed).\\

Exchangeability is possibly the most important concept in probability theory extending the notion of independence. Investigations in this direction were initiated in the 1930's simultaneously with the search for a general axiomatic approach to probability and started with a famous result of de Finetti (cf.\ \cite{DeFi1,DeFi2,Ki}) later extended by Hewitt and Savage \cite{HeSa}. A recent result in this area in connection with copulas is given by Mai and Scherer \cite{MaSch1}, where an interested reader may find an excellent overview of the subject together with extensive bibliography.

However, our aim is not so much to study exchangeability as the lack of it, a subject that had attracted little attention up to the point when Klement and Mesiar \cite{KM} and also Nelsen \cite{N} noticed it only a dozen of years ago using copulas.
In practice dependence is often asymmetric (i.e. non-exchangeable), as data collected from the real world may exhibit. This necessitates developing asymmetric copulas that can model such data and it also urges the study of various measures of asymmetry that may help the practitioners to decide about which copulas to choose in their models according to the data. Papers \cite{KM,N} started a vivid interest in the subject. De Baets, De Meyer, and Mesiar \cite{DeBDeMM} present an asymmetric version of semilinear copulas as an asymmetric version of the previously introduced symmetric version of semilinear copulas \cite{DKMS}. Durante, Klement, Sempi, and \'{U}beda-Flores \cite{DKSU-F} introduce a measure of asymmetry $\mu$ in general and $\mu_p$ for $p\in[1,\infty]$ in particular; we will recall these in Section \ref{preliminaries} of this paper.
It was shown in \cite{KM,N} that $\displaystyle \mu_\infty(C) \leq\frac13$ for any copula $C$ and that the bound is attained so that $\displaystyle\frac13$ is the sharp bound of asymmetry measure $\mu_\infty$. Similarly, the sharp bound of asymmetry measure $\mu_\infty$  for the set of positively quadrant dependent copulas was given in \cite{DeBDeMM} and the sharp bound of asymmetry measure for the set of negatively quadrant dependent copulas was given in \cite{DP1}. Some other techniques of constructing asymmetric copulas and finding sharp bounds of asymmetry measure for various families of copulas were presented in the papers mentioned above as well as in \cite{DP2}, \cite{DM}, and \cite{GN}.\\

One of the main contributions of this paper are exact expressions and systematic application of the \emph{maximal asymmetry function} on a given family of copulas. This function was first considered by Klement and Mesiar \cite{KM} and by Nelsen \cite{N} on the family of all copulas. It is computable for all the families that we are studying and it is helpful not only in determining the sharp bound of measures of asymmetry for a given family, but also in the detailed analysis of shocks in a given model at which the bound is attained. This brings us to another important contribution of ours, the stochastic interpretation of shocks in these models. We believe this approach is helpful to practitioners in search of the copulas that will fit their data the best and simultaneously the models they are hoping for. So, when one has measured the asymmetry of the data at hand and established a statistical hypothesis on the studied phenomenon, one can compare using our approach whether the pattern of occurrences of the shocks given the asymmetry allows for the desired conclusion.
\\

The paper is organized as follows: The main tools are presented in Sections \ref{preliminaries} and \ref{sec:pqd_nqd}. In Section \ref{sec:pqd_nqd} we use them to find the maximal asymmetry functions of the families of PQD and of NQD copulas. We also present families of copulas where the extremal values are attained. The analogous results on Marshall copulas are given in Section 4, the results on maxmin copulas in Section 5, and the results on RMM copulas in Section 6. Stochastic interpretations of the extremal families of copulas appearing in shock models are given in Section 7. All these results seem to be new. Stochastic interpretation of Marshall copulas of maximal asymmetry measure is illustrated in Figure \ref{slika11} (cf.\ also the comments just preceding it). Similarly, stochastic interpretations for maxmin copulas of maximal asymmetry measure are presented in Figures \ref{slika15}, \ref{slika16}, and for RMM copulas of maximal asymmetry measure in Figure \ref{slika23}, all the three explained in the paragraphs just preceding them. The development of these results turns out to be technically quite involved, so we postpone some of the proofs to the appendix.

\section{Preliminaries}\label{preliminaries}

A (bivariate) \emph{copula} is a function $C :[0, 1]^2 \rightarrow [0, 1]$ which satisfies:
\begin{description}
  \item[(a)] (\emph{boundary condition}) for every $x, y \in [0, 1], C(x, 0) = 0 = C(0, y)$ and $C(x, 1) = x$ and $C(1,y) = y$;
  \item[(b)] (\emph{2-increasingness}) for every $x_1, x_2, y_1, y_2 \in [0, 1]$ such that $x_1 \leq x_2$ and $y_1 \leq y_2$, $V_C(R) = C(x_2, y_2) - C(x_2, y_1) -$ $C(x_1, y_2)
	+ C(x_1, y_1) \geq 0$ (in other words, for all rectangles $R = [x_1, x_2]\times[y_1, y_2]$ whose vertices lie in $[0, 1]^2, V_C(R) \ge 0)$.
\end{description}
Sklar \cite{Sklar} showed that given any copula $C$ and arbitrary univariate distribution functions $F$ and $G$ function $H=C(F,G)$ is a joint distribution function, and that, even more importantly, given any joint distribution function $H$ of a random vector with marginal distribution functions $F$ and $G$ there is a copula $C$ such that $H=C(F,G)$. Note that for an exchangeable random vector $(X,Y)$ the components are, firstly, identically distributed so that $F=G$, and secondly, the exchange of the two components yields the same joint distribution, so that $C(F(x),F(y))= C(F(y),F(x))$ yielding that copula $C$ is \emph{symmetric} (i.e.\ $C(u,v)= C(v,u)$ for all $u,v\in [0,1]$). Conversely, given a symmetric copula $C$ and a univariate distribution function $F$, the joint distribution function $C(F,F)$ belongs to an exchangeable pair of random variables. So, from the copula point of view studying exchangeable random variables is equivalent to studying symmetric copulas.\\

Many classical copulas are symmetric, (sometimes also called \emph{exchangeable} due to reasons given above): Archimedean and meta-elliptical copulas are prime examples. Observe that two of the most important copulas, the \emph{Fr\'echet-Hoeffding lower bound}, respectively\emph{ upper bound}, $W(u,v)=\max\{0,u+v-1\}$, respectively $M(u,v)= \min\{u,v\}$ (being so called due to the fact that  $W(u,v) \leq C(u,v)\leq M(u,v)$ for every copula $C$ and all $u,v\in [0,1]$) are symmetric. Also, the independence of two random variables is being modeled via the \emph{product copula} $\Pi(u,v)=uv$ which is also symmetric. This is a copula argument showing the well known fact, and easy to see, that every independent pair of equally distributed random variables is exchangeable. So, in view of the classical exchangeability results one might vaguely think of more asymmetric copulas as modeling more dependent relations among random variables.\\

Let us briefly recall at this point some further notions to be needed in the sequel. Copula $C(u,v)$ is called \emph{positive quadrant dependent (PQD for short)} if $\Pi(u,v)\leq C(u,v)$ for all $u,v\in [0,1]$, and it is \emph{negative quadrant dependent (NQD for short)} if $C(u,v) \leq \Pi(u,v)$ for all $u,v\in [0,1]$. We denote by $\mathcal{C}$, respectively $\mathcal{P}$, respectively $\mathcal{N}$, the set of all copulas, respectively PQD copulas, respectively NQD copulas. For any $C\in\mC$ we denote by $C^t$ the copula defined by $C^t(u,v)=C(v,u)$ for all $u,v\in [0,1]$ (so that $C$ is symmetric if and only if $C=C^t$). The {maximal asymmetry function} for any particular family of copulas is defined as the point-wise supremum of all possible differences of $|C-C^t|$ when $C$ runs through the given family. Klement and Mesiar \cite{KM} were using this notion only on the family $\mathcal{C}$. The major role of this function in our paper has been described above. We also denote by $\widehat{C}$ the copula defined by $\widehat{C}(u,v) =u+v-1+C(1-u,1-v)$ for all $u,v\in [0,1]$, i.e.\ the \emph{survival copula} of $C$. We refer to monographs \cite{DS,Joe,MaSch2,Nelsen} for further details on copulas.

To quantify asymmetry the authors in \cite{DKSU-F} introduced the notion of a measure of asymmetry.
A function $\mu: \mC\to [0,\infty)$ is {\em a measure of asymmetry (or a measure of non-exchangeability) for copula $C$} if it
satisfies the following properties:
\begin{description}
  \item[(B1)] there exists $K\in [0,\infty)$ such that, for all $C\in\mC$ we have $\mu(C) \le K$,
\item[(B2)] $\mu(C) = 0$ if and only if $C$ is symmetric,
\item[(B3)] $\mu(C) = \mu(C^t)$ for every $C\in\mC$,
\item[(B4)] $\mu(C) = \mu(\widehat{C})$ for every $C\in\mC$,
\item[(B5)] if $(C_n)_{n\in\mathbb{N}}$ and $C$ are in $\mC$, and if $(C_n)_{n\in\mathbb{N}}$ converges uniformly to $C$, then $(\mu(C_n))_{n\in\mathbb{N}}$ converges to $\mu(C)$.
\end{description}

A large class of measures of asymmetry is provided in \cite[Theorem 1]{DKSU-F}: Let $d_p$ be the classical $L_p$ distance in $\mC$ for $p \in [1,\infty]$, viz., for all $A, B \in \mC$ we have
$$d_p(A, B) =\left(\int_0^1\int_0^1\left|A(u,v)-B(u,v)\right|^pdu\,dv\right)^{\frac1p},$$
when $p\in [1,\infty)$ and, for $p =\infty$, we have
$$d_{\infty}(A, B) = \max_{x,y\in [0,1]}\left|A(x,y)-B(x,y)\right|.$$
Then, for every $p \in [1,\infty]$, the measure of asymmetry $\mu_p: \mC\to [0,\infty)$ is given by
$$\mu_p(C) = d_p(C,C^t ).$$
%\medskip

\section{Positive and negative quadrant dependent copulas}\label{sec:pqd_nqd}

As pointed out in the introduction it is significant in applications to know what is the extremal value of asymmetry for a particular family of copulas $\mF\subset\mC$, i.e., what is the exact value of $\mu_p(\mF)= \sup_{C\in\mF}(\mu_p(C))$, or at least what is an upper bound for the supremum. The values of measures of asymmetry are of special interest for families of copulas characterized by some statistical or analytical properties. For $\mP$ and $\mN$ (defined in the introduction), say, it was shown in \cite{DeBDeMM}, respectively \cite{DP1}, that
\begin{equation}\label{muP_muN}
\mu_{\infty}(\mP)=3-2\sqrt{2},\quad\textrm{\ respectively,\ }\quad \mu_{\infty}(\mN)=\sqrt{5}-2,
\end{equation}
and either of the extreme values is attained.
\medskip

In addition to the measure of asymmetry, it is also of interest to know the maximal difference $\left|C(x,y)-C(y,x)\right|$ at a particular point $(x,y)\in [0,1]\times[0,1]$. Here we follow Klement and Mesiar \cite{KM} and define function $d^*_{\mF}:[0,1]\times[0,1]\to [0,\infty)$ that we call the \emph{maximal asymmetry function}. Its value at a fixed point $(x,y)\in [0,1]\times[0,1]$ is given by
$$d^*_{\mF}(x,y)=\sup_{C\in\mF}\left\{|C(x,y)-C(y,x)|\right\}.$$
%(We can write maximum instead of supremum since the right hand side is attained for each $(x,y)$.)
Klement and Mesiar showed that
\begin{equation}\label{d_star_C}
d^*_{\mC}(x,y)=\min\{x,y,1-x,1-y,|x-y|\}.
\end{equation}
and that the family of copulas $W_{\lambda}(x,y)=\max\{M(x,y-\lambda),W(x,y)\}$ for $\lambda\in [0,1]$ is such that $$d^*_{\mC}(x,y)=|W_{\lambda}(x,y)-W_{\lambda}(y,x)|$$
for all $x,y\in [0,1]$ with $|x-y|=\lambda$, see \cite[Proposition  2.2]{KM}.

They use this function only on the class $\mC$, while we will now compute the maximal asymmetry functions for the families $\mathcal{P}$ and  $\mathcal{N}$. As it turns out this is of independent interest, since we will later compare them with the maximal asymmetry functions of the classes of copulas arising in shock models.

Observe that equality (\ref{d_star_C}) tells us that the value of $d^*_{\mC}$ is equal to the the maximal possible difference, i.e., to $|M(x,y)-W(x,y)|$ when the point $(x,y)$ is far enough from the diagonal, while closer to the diagonal the value of $d^*_{\mC}$ is equal to $|x-y|$. Below we show that similar results hold also for the classes of PQD, respectively NQD copulas. Namely, the value of $d^*_{\mP}$ is equal to the the maximal possible difference, i.e., to $|M(x,y)-\Pi(x,y)|$ when the point $(x,y)$ is far enough from the diagonal, while closer to the diagonal the value of $d^*_{\mP}$ is equal to $|x-y|$. Similar result holds for $d^*_{\mN}$. Precise statements are given in Theorems \ref{d*PN1} and \ref{d*PN2} below.

\begin{theorem}\label{d*PN1}
Maximal asymmetry function of the family $\mathcal{P}$ of PQD copulas is equal to
$$d^*_{\mP}(x,y)=\min\{x(1-y),(1-x)y,|x-y|\}.$$

%Furthermore, the family of copulas
%$$P_{\lambda}(x,y)=\max\{M(x,y-\lambda),xy\}$$
%for $\lambda\in [0,1]$ is such that
%\begin{equation}\label{d_star_P}
%d^*_{\mP}(x,y)=|P_{\lambda}(x,y)-P_{\lambda}(y,x)|
%\end{equation}
%for all $x,y\in [0,1]$ with $|x-y|=\lambda$.
For $0\le a \le 1$ and $\frac{2a}{a+1}\le b \le 1$ the bound $d^*_{\mP}(a,b)$ is attained, for example, by PQD copula
$$P_{a,b}(x,y)=\left\{\begin{array}{ll}
        x; & \text{if } 0\le x\le 1-b+a \text{ and }x+b-a\le y\le 1, \\
        xy+(1-b+a-x)(y-b+a); & \text{if } 0\le x\le 1-b+a \text{ and }b-a\le y\le x+b-a, \\
        xy; & \text{ otherwise. }
        \end{array}\right. $$
For  $0< a < 1$ and $a\le b \le \frac{2a}{a+1}$ the bound $d^*_{\mP}(a,b)$ is attained, for instance, by absolutely continuous PQD copula with density
\begin{equation}\label{banana}
q_{a,b}(x,y)=\left\{\begin{array}{ll}
        \frac{2a-b}{a^2}; & \text{if } 0\le x,y\le a, \\
        \frac{1}{a}; & \text{if } 0\le x\le a \text{ and }a< y\le b, \\
        \frac{1}{1-b}; & \text{if } a < x\le b \text{ and }b < y\le 1, \\
        \frac{b-a}{a(1-b)}; & \text{if } b < x\le 1 \text{ and }0\le y\le a, \\
        \frac{1+a-2b}{(1-b)^2}; & \text{if } b < x,y\le 1, \\
        0; & \text{ otherwise. }
        \end{array}\right.
\end{equation}
\end{theorem}

\prf
Choose $C\in\mP$ and assume that $y\ge x$. Then we have $xy\le C(x,y)\le x$. It follows that $|C(x,y)-C(y,x)|\le x-xy=x(1-y).$ Similarly, we see that $|C(x,y)-C(y,x)|\le y-xy=y(1-x)$ if $x\ge y$. We also have $d^*_{\mP}(x,y)\le d^*_{\mC}(x,y)$ for all $x$ and $y$. Relation (\ref{d_star_C}) implies that $d^*_{\mP}(x,y)\le d^*_{\mC}(x,y)\le |x-y|$, so it follows that
\begin{equation}\label{ineq}
d^*_{\mP}(x,y)\le \min\{x(1-y),(1-x)y,|x-y|\}.
\end{equation}
We complete the proof by showing that the inequality in (\ref{ineq}) is attained for each point in $[0,1]^2$.

Assume that $(a,b)\in [0,1]^2$ and that $a\le b$. The case $a\ge b$ follows by symmetry. For $\frac{2a}{a+1}\le b\le 1$ observe that copulas $P_{a,b}$ are PQD and that $P_{a,b}(a,b)-P_{a,b}(b,a)=a-ab=d^*_{\mP}(a,b)$.

For $a\le b \le \frac{2a}{a+1}$ we denote the copula with density $q_{a,b}(x,y)$ by $Q_{a,b}$. Under the assumed conditions on $a$ and $b$ one can check that all the expressions in $a$ and $b$ on the right-hand side of (\ref{banana}) are nonnegative. Observe that $Q_{a,b}(a,a)=2a-b\ge \frac{2a^2}{a+1}\ge a^2$. Since $a-b^2\ge\left(\frac{1-a}{1+a}\right)^2$ it follows that $Q_{a,b}(b,b)=a\ge b^2$. These, together with the fact that the density is uniform on each of the rectangles where it is nonzero, implies that $Q_{a,b}$ is PQD. Finally, it is straightforward to check that $Q_{a,b}(a,b)-Q_{a,b}(b,a)=a-(2a-b)=b-a=d^*_{\mP}(a,b)$.
\qed\\

\begin{figure}[h]
\begin{center}
\includegraphics[width=6cm]{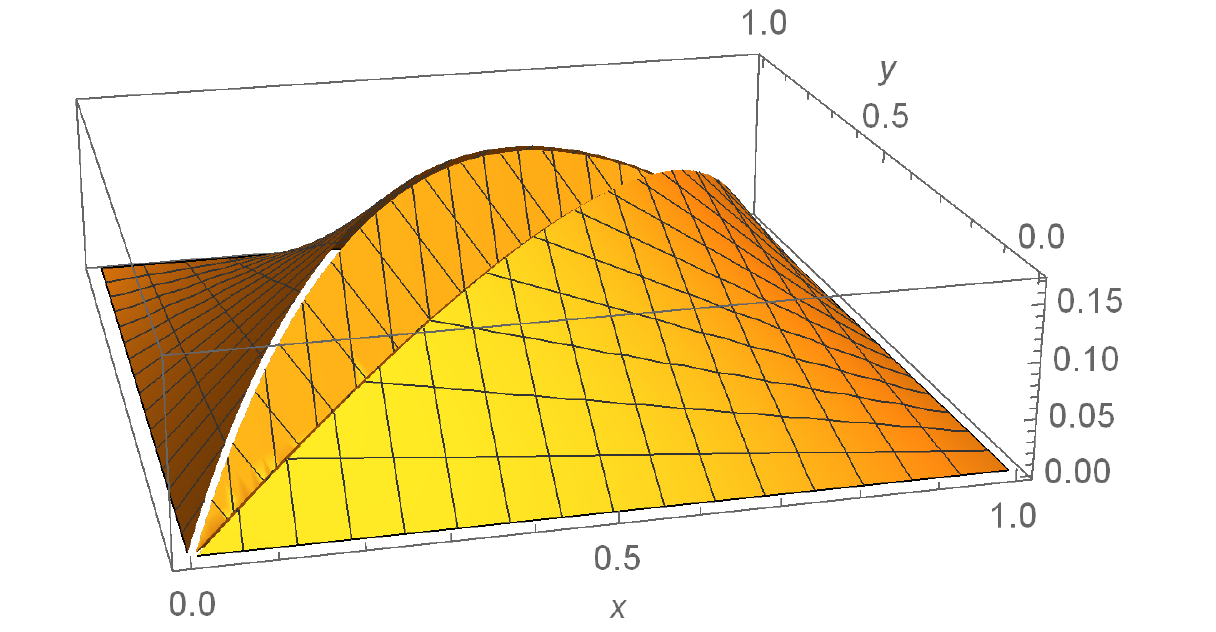} \quad \quad
\includegraphics[width=6cm]{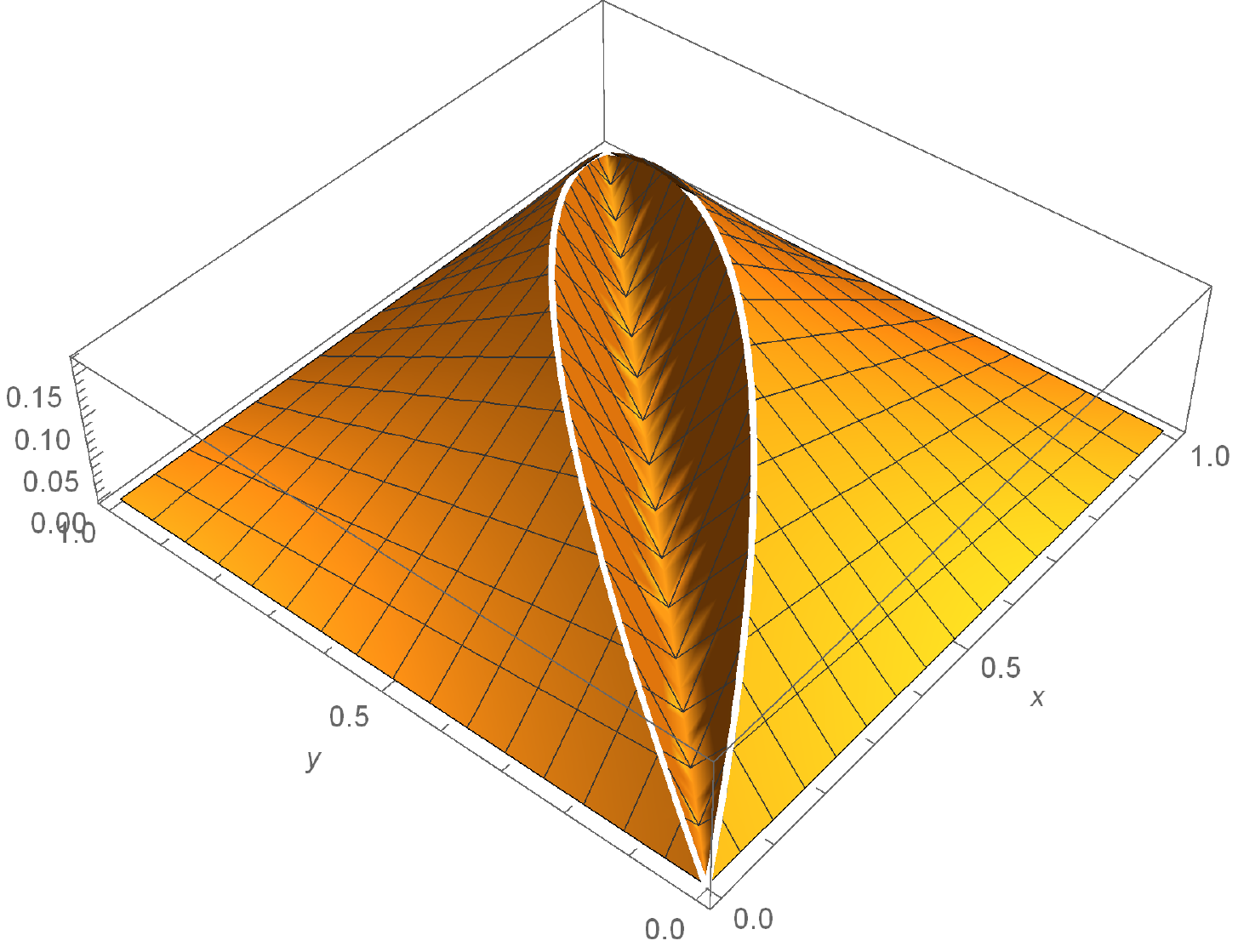}
\caption{Graph of $d^*_{\mP}$.} \label{slika01}
\end{center}
\end{figure}

Figure \ref{slika01} shows two views of the maximal asymmetry function for $\mP$.

%Observe that the family of copulas $P_{\lambda}$, where the maximal values of $d^*_{\mP}$ are attained, is %obtained from the family of copulas $W_{\lambda}$ by replacing the lower bound $W$ for the family of all %copulas $\mC$ with the lower bound $\Pi$ for the family of all PQD copulas $\mP$.

\begin{theorem}\label{d*PN2}
Maximal asymmetry function of the family  $\mathcal{N}$ of NQD copulas is equal to
$$d^*_{\mN}(x,y)=\min\{xy,(1-x)(1-y),|x-y|\}.$$
Furthermore, the family of copulas
$$N_{\lambda}(x,y)=\max\{W(x,y),\min\{y-\lambda,xy\}\}$$
for $\lambda\in \left[0,\frac12\left(3-\sqrt{5}\right)\right]$ is such that
\begin{equation}\label{d_star_N}
d^*_{\mN}(x,y)=|N_{\lambda}(x,y)-N_{\lambda}(y,x)|
\end{equation}
for all $x,y\in [0,1]$ with $|x+y-1|=\mu$, where $\mu=\frac{1-3\lambda+\lambda^2}{1-\lambda}$ and $\mu\in [0,1]$.
\end{theorem}

\prf
Consider the class of NQD copulas and choose $C\in\mN$. We assume that $x+y\le 1$. Then we have $0\le C(x,y)\le xy$ and so $|C(x,y)-C(y,x)|\le xy.$ Now, assume that $x+y\ge 1$. Then $x+y-1\le C(x,y)\le xy$ and thus $|C(x,y)-C(y,x)|\le xy-x-y+1=(1-x)(1-y).$ We also have $d^*_{\mN}(x,y)\le d^*_{\mC}(x,y)$ for all $x$ and $y$. Hence (\ref{d_star_C}) implies
$$ d^*_{\mN}(x,y)\le \min\{xy,(1-x)(1-y),|x-y|\}.$$

Observe that for $\lambda\in \left[0,\frac12\left(3-\sqrt{5}\right)\right]$ we have $\mu\in [0,1]$, and that
$$N_{\lambda}(x,y)=\left\{\begin{array}{ll}
        xy; & \text{if } 0\le x\le 1-\lambda \text{ and }\frac{\lambda}{1-x}\le y\le 1, \\
        y-\lambda; & \text{if } 0\le x\le 1-\lambda \text{ and }\lambda\le y\le \frac{\lambda}{1-x}, \\
        x+y-1; & \text{if } 1-\lambda\le x\le 1 \text{ and } 1-x\le y\le 1,\\
        0; & \text{ otherwise. }
        \end{array}\right. $$
A straightforward calculation shows that for $\lambda\in \left[0,\frac12\left(3-\sqrt{5}\right)\right]$, $x\in [0,1-\mu]$ and $y=1-\mu-x$ we have  $|N_{\lambda}(x,y)-N_{\lambda}(y,x)|=\min\{xy, |y-x|\}.$ Furthermore, equality  $|N_{\lambda}(x,y)-N_{\lambda}(y,x)|=\min\{(1-x)(1-y), |y-x|\}$ holds for $x\in [1-\mu,1]$ and $y=1+\mu-x$. Hence (\ref{d_star_N}) follows.
\qed\\

\begin{figure}[h]
\begin{center}
\includegraphics[width=6cm]{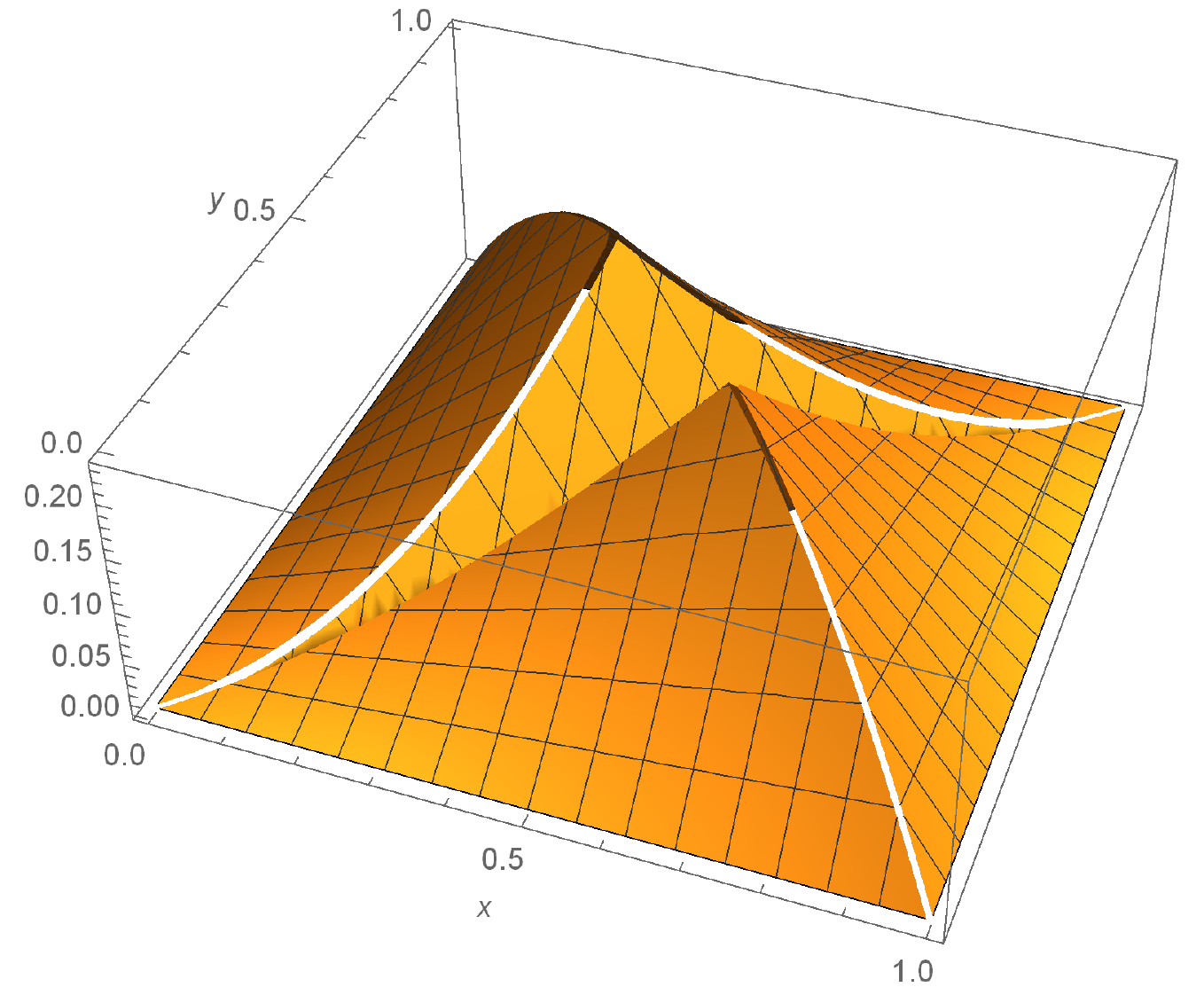} \quad \quad
\includegraphics[width=6cm]{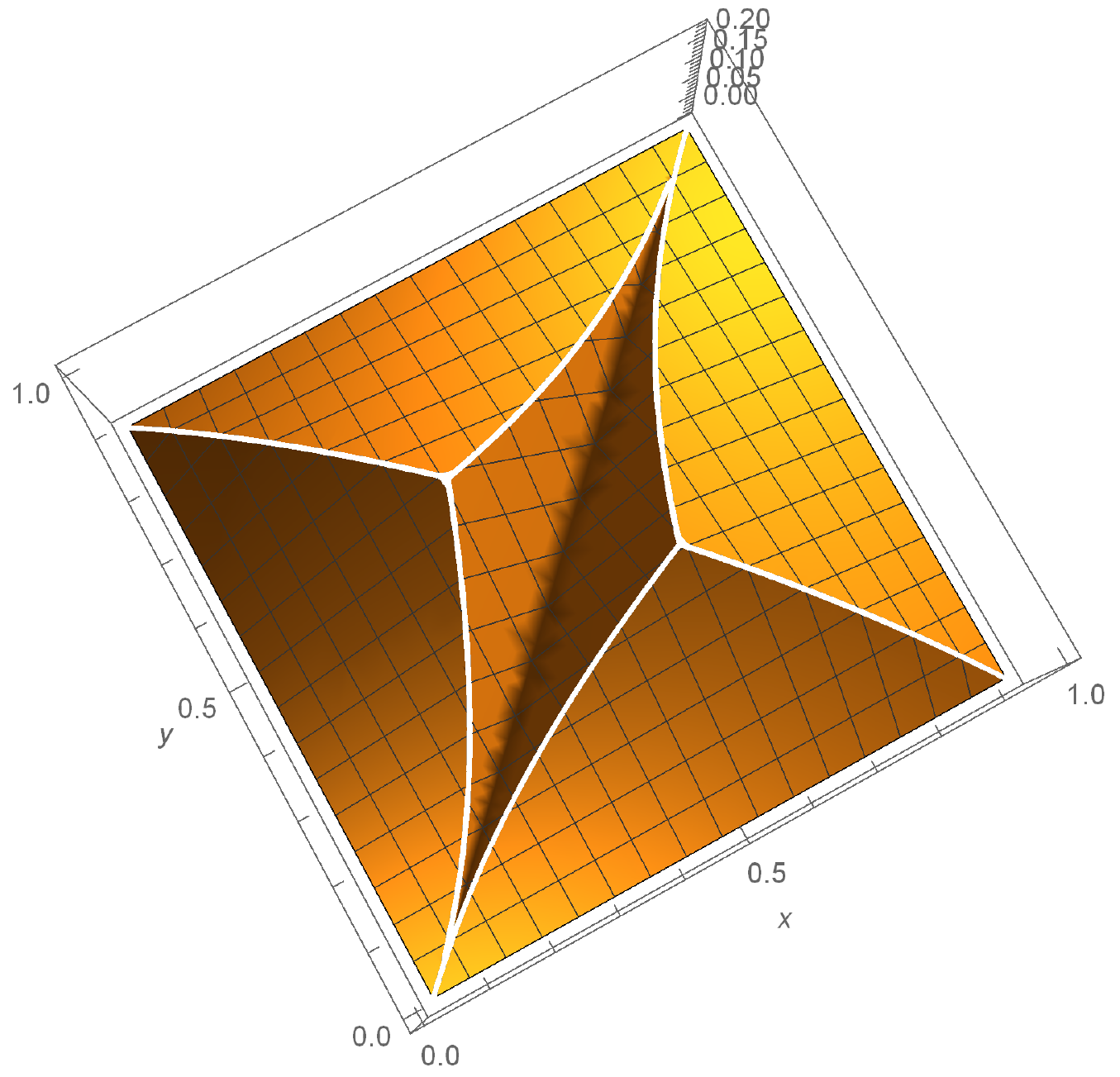}
\caption{Graph of $d^*_{\mN}$.} \label{slika03}
\end{center}
\end{figure}

Figure \ref{slika03} shows two views of the maximal asymmetry function for $\mN$.

Observe that the family of copulas $N_{\lambda}$, where the maximal values of $d^*_{\mN}$ are attained, is obtained from the family of copulas $W_{\lambda}$ using the upper bound $\Pi$ for the family of all NQD copulas $\mN$. Namely, we have $N_{\lambda}(x,y)=\max\{W(x,y),\min\{M(x,y-\lambda),xy\}\}=\max\{W(x,y),\min\{y-\lambda,xy\}\}$, since $xy\le x$ for $y\in [0,1]$.

\begin{remark}
The extremal values (\ref{muP_muN}) for $\mP$ and $\mN$, respectively, are attained at $(x,y)=(\sqrt{2}-1,2-\sqrt{2})=:(a_0,b_0)$ and $\left(\frac{3-\sqrt{5}}{2},\frac{\sqrt{5}-1}{2}\right)$, respectively (see \cite{DeBDeMM,DP1}). %There, we have $$|x-y|=3-2\sqrt{2}=:\lambda_0\ \text{and}\ |x+y-1|=0=:\mu_0,$$ respectively.
Thus, the members from the families given in Theorems \ref{d*PN1} and \ref{d*PN2}, where the maximal values of $d^*_{\mP}$ and $d^*_{\mN}$, respectively, are attained, are $P_{a_0,b_0}$, $Q_{a_0,b_0}\in\mP$ and $N_{\lambda_0}\in\mN$, respectively. Here we have $$\frac{3-\sqrt{5}}{2}+\frac{\sqrt{5}-1}{2}-1=0=:\mu_0\text{ and }\lambda_0=\frac12\left(3-\mu_0-\sqrt{\mu^2_0-2\mu_0+5}\right)=\frac{3-\sqrt{5}}{2}.$$
\qed
\end{remark}

\section{ Marshall copulas }\label{Marshall_copulas}

In this section we study the asymmetry of Marshall copulas. Before we recall their definition we introduce some notation. For a nondecreasing function $f: [0, 1] \to [0, 1]$ we write $$f^*(x)=\frac{f(x)}{x},$$
where $f^*(0) =\lim_{x\downarrow 0}\frac{f(x)}{x}$ if the limit exists or $f^*(0) =\infty$ otherwise. \\[1mm]

A \emph{Marshall copula} is a copula given by
$$C(x, y) = \min\{yf(x), xg(y)\} ,$$
where $f, g: [0, 1] \to [0, 1]$ are nondecreasing functions that satisfy the following properties:
\begin{description}
  \item[(a)] $f(0) = g(0) = 0$, $f(1) = g(1) = 1$,
  \item[(b)] $f^*$ and $g^*$ are nonincreasing.
\end{description}
See \cite{M} and \cite{DKMS}.

Note that the conditions imposed on generators $f$ and $g$ of a Marshall copula imply that they are continuous functions everywhere on $(0,1]$ that are possibly discontinuous at $0$. For instance, copula $M$ is a symmetric Marshall copula (i.e. a semilinear copula by \cite{DKMS}) with generator
$$f(t)=g(t)=\left\{\begin{array}{ll}
        0; & \text{if } x=0, \\
        1; & \text{if } x>0.
        \end{array}\right. $$
\medskip

Since the expression $|C(x, y) - C(y, x)|$ is symmetric with respect to the line $x=y$, it suffices to study only the case $y \ge x$ { in all considerations that follow}.

The proof of the following result is a bit technical. We include it in the Appendix.

\begin{lemma} \label{ocenaM}\label{zgornja_mejaM}
Let $C$ be a Marshall copula and $x, y \in [0, 1]$ with $y \ge x$. Then
$$|C(x, y) - C(y, x)| \le F(x,y)$$
where
$$F(x, y) = \left\{\begin{array}{ll}
        x(1-y); & \text{if } x \le y^2, \\
        \frac{x(y-x)}{y}; & \text{if } x > y^2. \\
        \end{array}\right.$$
In, particular, we have that
$$F(x, y) \le \frac{4}{27}$$
for all $(x,y)$ and
$$F\left(\frac49,\frac23\right) =F\left(\frac23,\frac49\right) = \frac{4}{27}.$$
\end{lemma}

\medskip

\begin{figure}[h]
\begin{center}
\includegraphics[width=7cm]{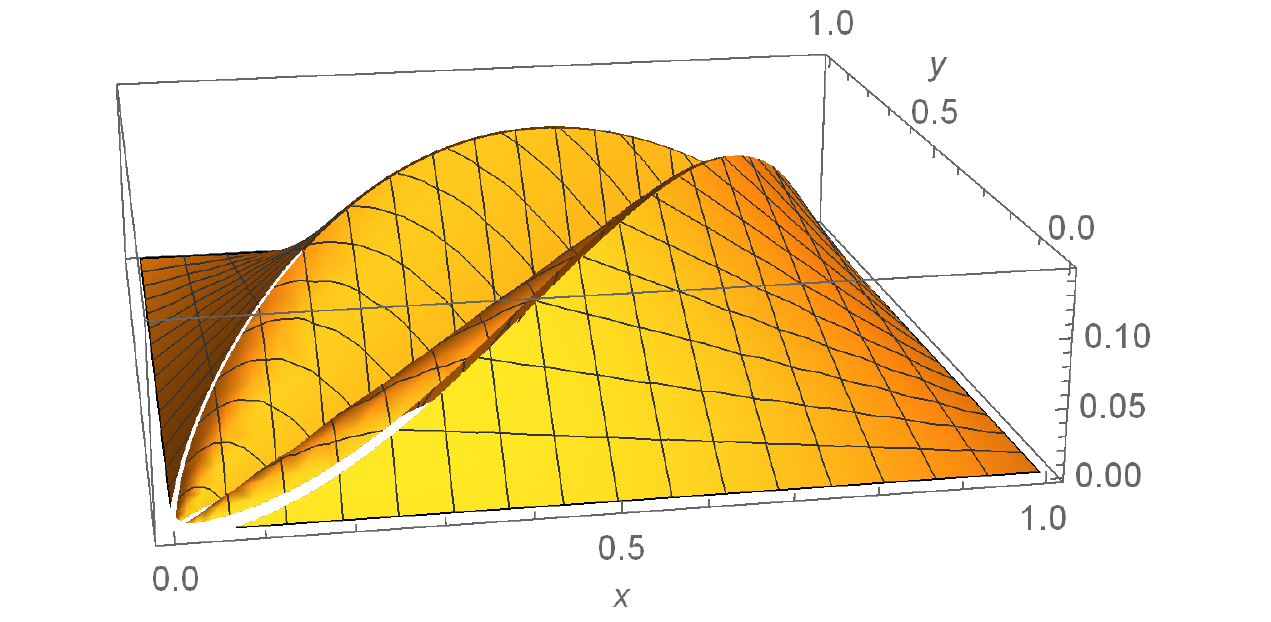} \quad \quad
\includegraphics[width=7cm]{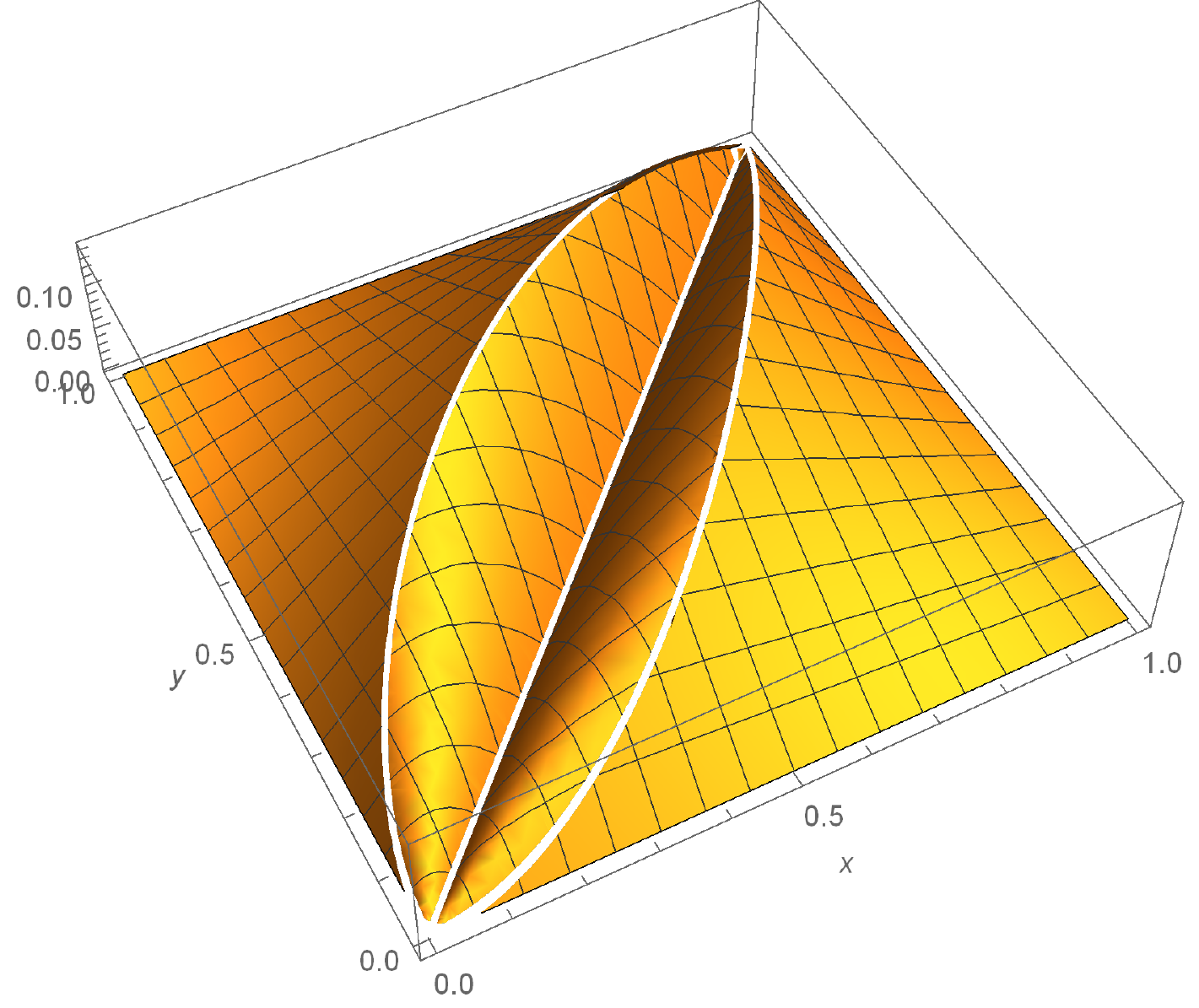}
\caption{Graph of $d^*_{\mM}$ for the set of Marshall copulas $\mM$.}\label{slika05}
\end{center}
\end{figure}

Figure \ref{slika05} shows two views of the upper bound for the expression $|C(x, y) - C(y, x)|$ for any Marshall copula $C$.  As we will see this is the maximal asymmetry function for the set of Marshall copulas.

\begin{example} \label{dosezenoM}
We consider functions
$$ f_{\mu}(x) = \left\{\begin{array}{ll}
        0; & \text{if } x = 0,\\
        \mu; & \text{if } 0 < x \le \mu, \\
        x; & \text{if } x > \mu,
        \end{array}\right.
\ \ \text{and} \ \
g(x) = \left\{\begin{array}{ll}
        0; & \text{if } x=0, \\
        1; & \text{if } x>0,
        \end{array}\right. $$
where $\mu\in (0,1)$.
Obviously $f_{\mu}$ and $g$ are nondecreasing. It is also easy to see that functions
$$ f_{\mu}^*(x) = \left\{\begin{array}{ll}
        \infty; & \text{if } x = 0,\\
        \frac{\mu}{x}; & \text{if } 0 < x \le \mu, \\
        1; & \text{if } x > \mu,
        \end{array}\right.
\ \ \text{and} \ \
g^*(x) = \left\{\begin{array}{ll}
        \infty; & \text{if } x=0, \\
        \frac{1}{x}; & \text{if } x>0,
        \end{array}\right.$$
are nonincreasing. So, $C_{\mu}(x,y)=\min\{yf_{\mu}(x),xg(y)\}$ is a Marshall copula. It is given by
\begin{equation}\label{C_mu}
    C_{\mu}(x, y) = \left\{\begin{array}{ll}
        x; & \text{if } 0\le {x}\le\mu y, \\
        \mu y; & \text{if } \mu{y} \le x \le \mu, \\
        xy; & \text{if } \mu \le x \le 1.
        \end{array}\right.
\end{equation}
For $x\le y$, we have
\begin{equation}\label{DeltaC_mu}
|C_{\mu}(x, y) - C_{\mu}(y,x)|= \left\{\begin{array}{ll}
        x(1-\mu); & \text{if } 0\le {x}\le\mu y,0\le y\le \mu, \\
        x(1-y); & \text{if } 0\le {x}\le\mu y,\mu\le y\le 1, \\
        \mu (y-x); & \text{if } \mu{y} \le x \le y, 0\le y\le \mu, \\
        y(\mu-x); & \text{if } \mu y \le x\le\mu, \mu\le y \le 1, \\
        0; & \text{if } \mu\le x\le y, \mu\le y\le 1.
        \end{array}\right.
\end{equation}
Thus, we have the equality $|C_{\mu}(x, y) - C_{\mu}(y,x)|=F(x,y)$ in the region $R_{\mu}=\{(x,y);\ \mu\le y\le 1, 0\le {x}\le\mu y\}$. Also, we have $\mu(y-x)=\frac{x}{y}(y-x)$ if and only if $y=\frac{x}{\mu}$. So, the difference  $|C_{\mu}(x,y)-C_{\mu}(y,x)|$ and $F(x,y)$ coincide on the segment $I_{\mu}=\{(x,\frac{x}{\mu});\ 0\le x\le \mu^2\}$, too.
%if $(x,y)\in [0,\mu^2]\times [\mu,1]$. % or  $(x,y)\in [\mu,1]\times [0,\mu^2]$.
For $\mu=\frac23$, we obtain
$$\mu_\infty(C_{\frac{2}{3}}) = \frac{4}{27}.$$

Let us observe in passing that copulas $C_{\mu}$ are ordinal sums of copulas $M$ and $\Pi$, and so they are  special cases of more general patchwork construction of copulas. We refer to \cite[pp. 63-64]{Nelsen} for the notion of an ordinal sum of copulas and to \cite{DS-PS} for the general patchwork construction for bivariate copulas.
\qed
\end{example}

\medskip

%Figure \ref{slika07} shows the graph of a copula $C_{\frac23}$, its scatterplot of 2000 points, and the graph of the expression $|C_{\frac23}(x, y) - C_{\frac23}(y, x)|$, respectively:

\begin{figure}[h]
\begin{center}
\includegraphics[width=5cm]{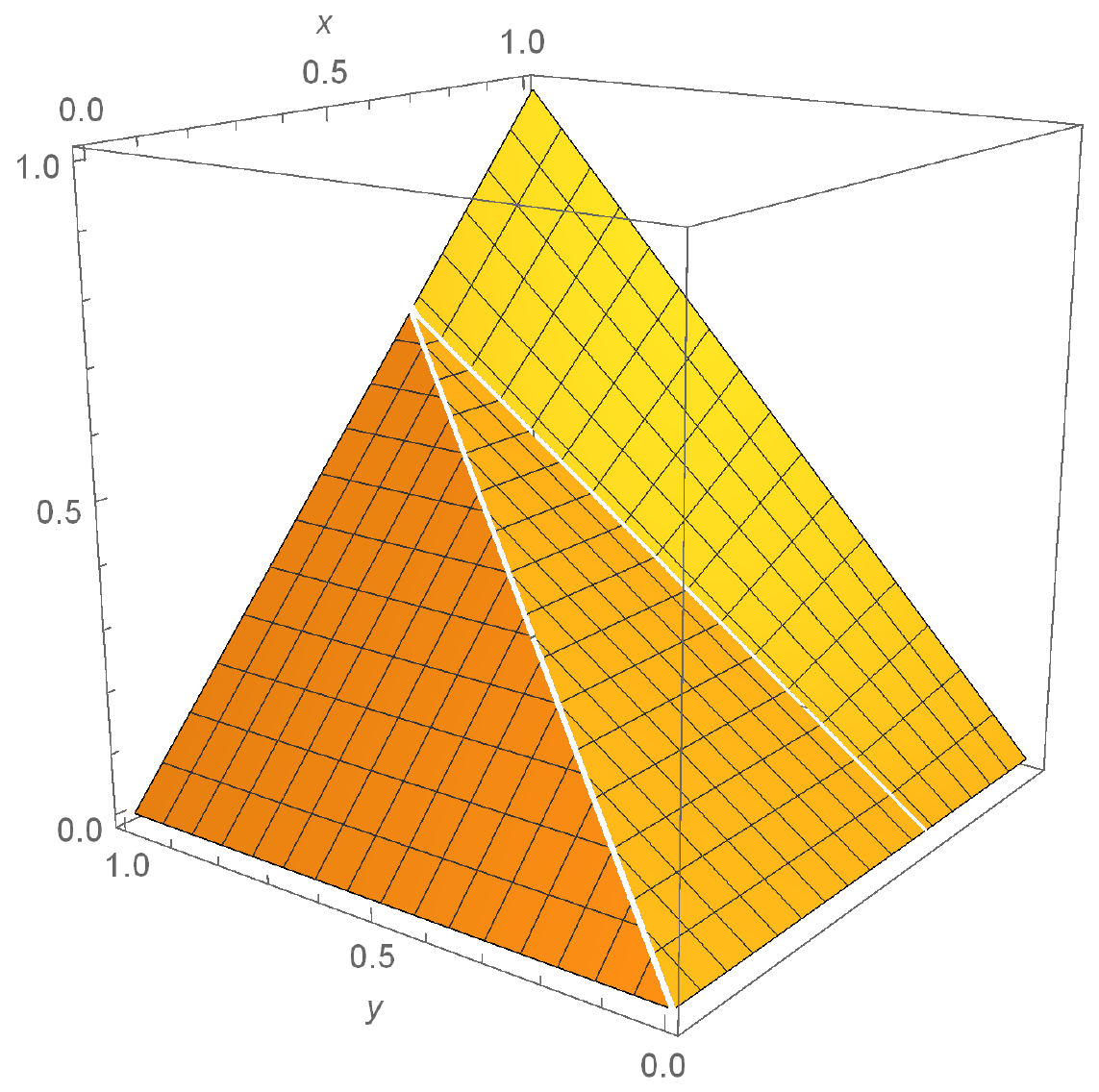} \ \
\includegraphics[width=42mm]{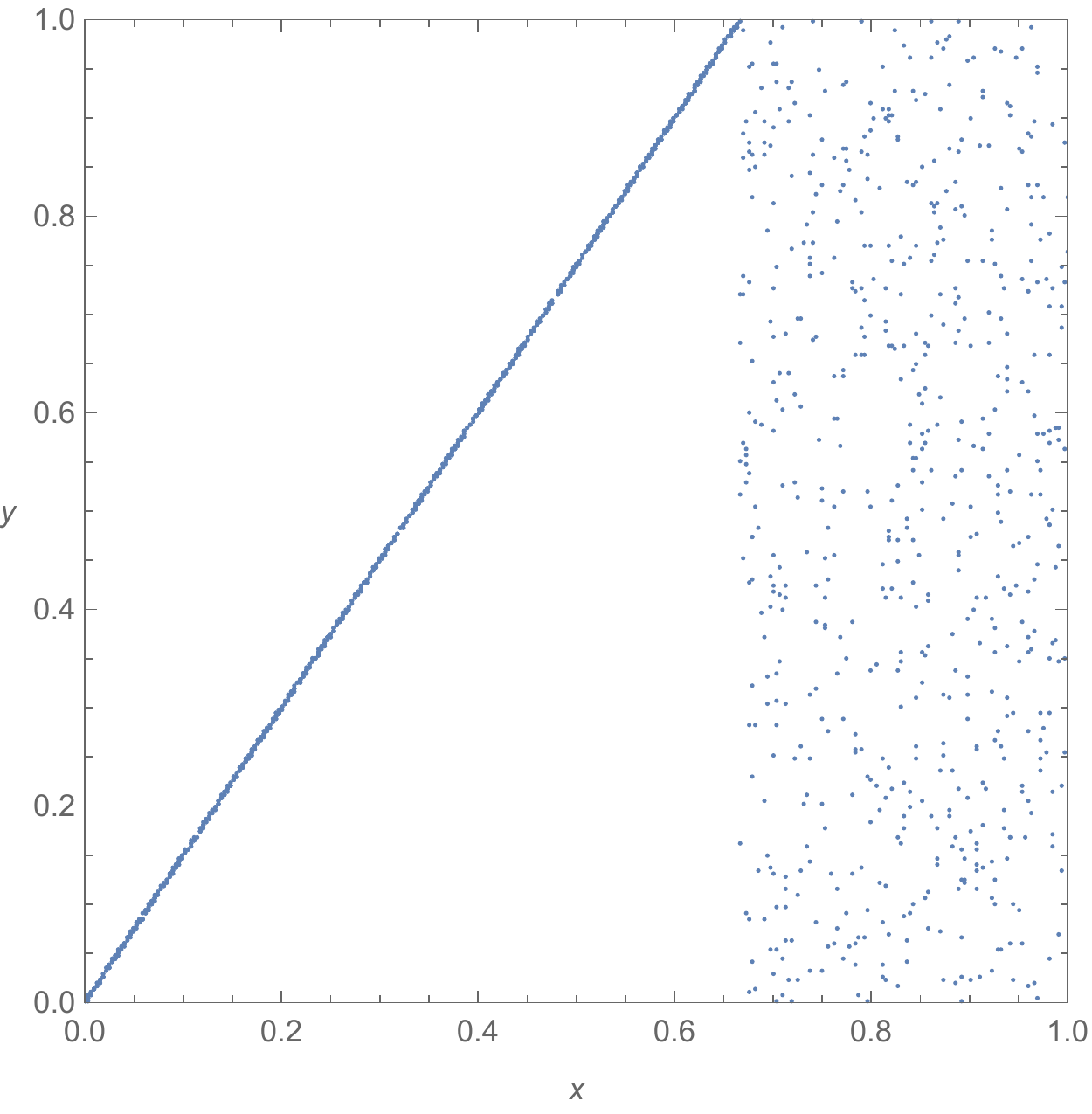} \ \
\includegraphics[width=70mm]{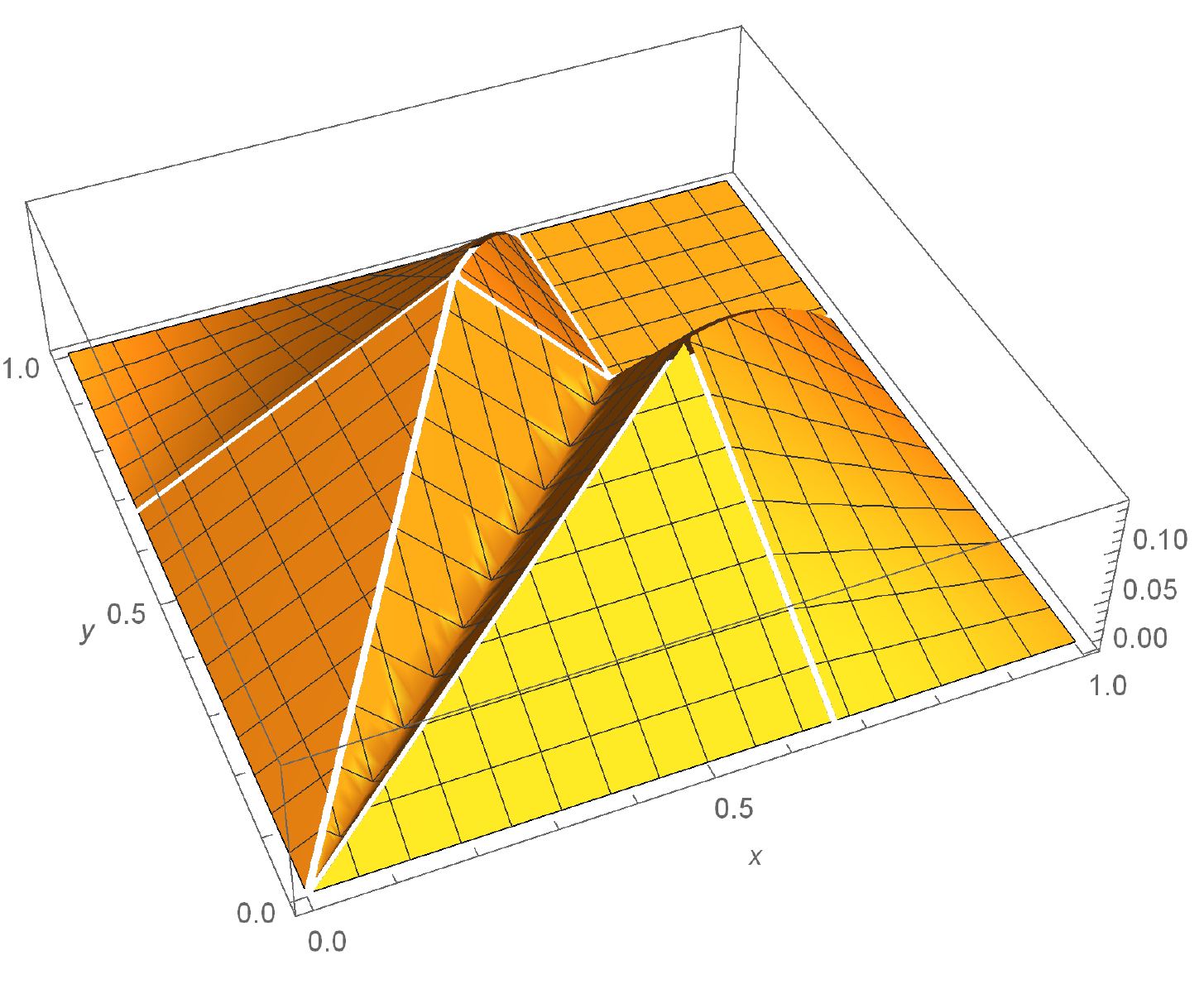}
\caption{Figure shows the graph of a copula $C_{\frac23}$, its scatterplot of 2000 points, and the graph of the expression $|C_{\frac23}(x, y) - C_{\frac23}(y, x)|$, respectively.}\label{slika07}
\end{center}
\end{figure}

The following is our main result for Marshall copulas. It is an immediate consequence of the preceding results of this section.

\begin{theorem}\label{d* for Marshall}
Let $\mM$ be the set of all Marshall copulas and $C\in \mM$. Then $\mu_\infty(C) \le \frac{4}{27}$ and the bound is attained, i.e., $\mu_{\infty}(\mM)=\frac{4}{27}$. The maximal asymmetry function for $\mM$ is equal to
$$d^*_{\mM}(x,y)= \left\{\begin{array}{ll}
        x(1-y); & \text{if } x \le y^2,\ \\
        \frac{x}{y}(y-x); & \text{if }  y^2\le  x\le y,\\
        y(1-x); & \text{if } y \le x^2,\ \\
        \frac{y}{x}(x-y); & \text{if }  x^2 \le y\le x.
        \end{array}\right.$$
\end{theorem}

\prf
The first claim follows from Lemma \ref{zgornja_mejaM} and Example \ref{dosezenoM}. The second claim follows by Example \ref{dosezenoM}, since the union of all the regions $R_{\mu}$ for $\mu\in (0,1)$ is the region $\{(x,y);\ 0\le y\le 1,\, 0\le x\le y^2\}$ and the union of all the segments $I_{\mu}$ for $\mu\in (0,1)$ is the region $\{(x,y);\ 0\le y\le 1,\, y^2\le x\le y\}$. We exchange the roles of $x$ and $y$ if $x\ge y$.
\qed

\begin{remark}
Since Marshall copulas are PQD we have $d^*_{\mM}(x,y)\le d^*_{\mP}(x,y)$ for all $(x,y)$. Observe that the two maximal asymmetry functions that are given in Theorems \ref{d*PN1} and \ref{d* for Marshall}, respectively, coincide for $(x,y)$ with $x\le y^2$ or $x\ge\sqrt{y}$.
\end{remark}

\begin{observation} \label{zg_mejaM_p}
Let $C$ be a Marshall copula and $p\in[1,\infty)$.
Then
\begin{equation}\label{Mp_ocena}
\mu_p(C) \le \left(\frac{2(2p + 3) B(p + 1, 2p + 3)}{p^2 + 3p + 2}\right)^{1/p}.
\end{equation}
In particular, $\mu_1(C) \le \frac{1}{18} \approx 0.0556$ and $\mu_2(C) \le \frac{1}{6\sqrt{6}} \approx 0.0680$.

If $C$ is a Marshall copula it follows from  Lemma \ref{ocenaM} that
$$ \mu_p(C) \le \left(2\int\!\!\int_T F(x,z)^p\,dx\,dz\right)^{1/p} $$
where $T$ is triangle $T = \{(x, y);\ 0\le x\le y\le 1\}$.
Using the Mathematica software \cite{Mathematica} we verified that the value of the integral is equal to the bound in (\ref{Mp_ocena}).
\end{observation}

\medskip

Recall that incomplete beta function is defined as
$$B(z; \alpha, \beta) = \int_0^z x^{\alpha-1}(1-x)^{\beta-1}\,dx.$$

\medskip

\begin{figure}[h]
\begin{center}
\includegraphics[width=9cm]{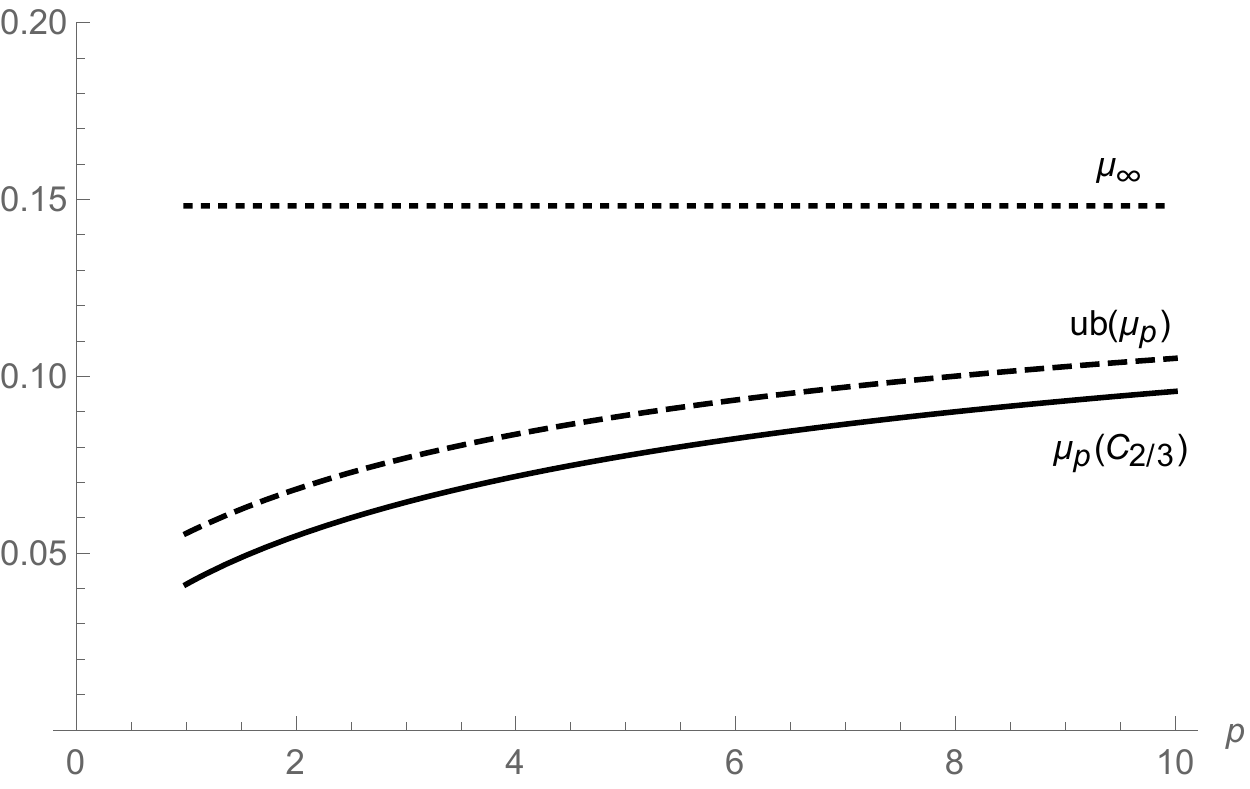}
\caption{Bounds for $\mu_p(\mM)$.}\label{slika10}
\end{center}
\end{figure}

\medskip

\begin{example} \label{dosezenoM_p}
Let $C_{\frac23}$ be  the Marshall copula from Example \ref{dosezenoM} { at $\mu=\frac23$} and $p \ge 1$.
Then $$\mu_p(C_{\frac23}) = \frac{2}{3}\left(\frac{8}{3(p+1)}B(\textstyle\frac{1}{3}; p+2, p+1) + \displaystyle\frac{2^{p+3}(4p+5)}{3^{2p+3}(p+1)^2(p+2)}\right)^{1/p}.$$
In particular, $\mu_1(C_{\frac23}) = \frac{10}{243} \approx 0.0411$ and $\mu_2(C_{\frac23}) = \frac{2}{81}\sqrt{\frac{74}{15}} \approx 0.0548$.

We integrated the function $|C_{\frac23}(x,y)-C_{\frac23}(y,x)|$ over the triangle $T$ using the Mathematica software \cite{Mathematica} and obtained the value given above.

Thus, we have found an upper bound for asymmetry $\mu_p(\mM)$ for the set $\mM$ of Marshall copulas,
and $\mu_p(C_{\frac23})$ for the copula $C_{\frac23}$. They are not the same, but they are not far apart.

Figure \ref{slika10} shows the graphs of both as functions of $p$ when $1 \le p \le 10$ (full curve and dashed curve), and also the asymptote of both, $y = \frac{4}{27}$ (dotted line).
\end{example}

\section{ Maxmin copulas }

A \emph{maxmin copula} is a copula given by
\begin{equation}\label{maxmin}
C(x, y) = \min\{x, \phi(x)y - \phi(x)\psi(y) + x\psi(y)\},
\end{equation}
where $\phi, \psi: [0, 1] \to [0,1]$ are nondecreasing functions
that satisfy the properties:
\begin{description}
  \item[(a)] $\phi(0) = \psi(0) = 0$, $\phi(1) = \psi(1) = 1$,
  \item[(b)] $\phi^*(x) = \frac{\phi(x)}{x}$ and $\psi_* = \frac{1-\psi(x)}{x-\psi(x)}$
are nonincreasing.
\end{description}

In (b) above we have $\psi_*(x) = \infty$ if $\psi(x) = x$ for $x < 1$ and $\psi_*(1)=1$.

For further details on maxmin copulas confer \cite{OR} and \cite{DOOR}.

\medskip

Next we define another pair of auxiliary generating functions that will have an important role in the remainder of the paper. Let
\begin{equation}\label{f and g}
f(x) = \phi(x) - x \ \ {\rm and} \ \  g(x) = 1 - x - \psi(1 - x).
\end{equation}
Note that $f$ and $g$ are defined on $[0,1]$.
We denote by $\mF$ the set of functions $f: [0, 1] \to [0, 1]$ that satisfy the properties:
\begin{description}
  \item[(a)] $f(0) = f(1) = 0$,
  \item[(b)] $f(x) + x$  is nondecreasing,
  \item[(c)] $f^*(x)=\frac{f(x)}{x}: (0, 1] \to [0, \infty)$ is nonincreasing.
\end{description}

\medskip

The following result is proved in \cite[Theorem  2]{KO}.

\begin{lemma}\label{takojzaF}
$C$ is a maxmin copula generated by functions $\phi$ and $\psi$ if and only if functions $f$ and $g$ defined by \eqref{f and g} belong to the set ${\mathcal F}$.  In this case we have $C(x,y)=\min\{x, xy+f(x)g(1-y)\}.$
\end{lemma}
\medskip

If $\phi$ and $\psi$ are related to $f$ and $g$ via (\ref{f and g}) then we say either that $\phi$ and $\psi$ generate maxmin copula $C$ of (\ref{maxmin}) or that $f$ and $g$ generate $C$.
\medskip

\begin{lemma} \label{C-f-g}
Let $C$ be a maxmin copula generated by $f$ and $g$ from $\mF$. Given $x, y \in [0, 1]$ write $z = 1-y$. Then
\begin{equation}\label{deltamaxmin}
|C(x, y) - C(y, x)| = |\min\{xz, f(x)g(z)\} - \min\{(1-x)(1-z), f(1-z)g(1-x)\}|.
\end{equation}
\end{lemma}

\prf
Suppose $C$ is given by (\ref{maxmin}) and $f$, $g$ by (\ref{f and g}). Then a direct calculation yields
\begin{eqnarray*}
|C(x,y)-C(y,x)|&=&|\min\{x,xy+f(x)g(1-y)\}-\min\{y,xy+f(y)g(1-x)\}|\\
&=&|\min\{ x z,f(x)g(z)\}-\min\{(1-x)(1-z),f(1-z)g(1-x)\}|.
\end{eqnarray*}
\qed

\medskip

\begin{remark}\label{symmetry_mm}
Note that it follows from Lemma \ref{C-f-g} that the expression $|C(x, y) - C(y, x)|$ is symmetric with respect to the line $x=z$, i.e., $x+y=1$. Evidently, it is symmetric also with respect to the line $x=y$.

Due to the latter symmetry we will restrict our consideration only to the case $y\ge x$. Then we have that $x+z\le 1$.\qed
\end{remark}
\medskip

The next result gives an upper bound for expression (\ref{deltamaxmin}). Its proof is rather technical and it is included in Appendix.

\begin{lemma} \label{ocenaMM}
Let $x, z \in [0, 1]$ and $x+z \le 1$. For $f,g\in\mF$ we have
$$|\min\{xz, f(x)g(z)\} - \min\{(1-x)(1-z), f(1-z)g(1-x)\}| \le G(x,z)$$
where
$$G(x, z) = \left\{\begin{array}{ll}
        \frac{x(1-x-z)}{1-z}; & \text{if } 1 - \sqrt{x} \le z \le \min\{x, 1-x\}, \\
        \frac{z(1-x-z)}{1-x}; & \text{if } 1 - \sqrt{z} \le x \le \min\{z, 1-z\}, \\
        xz; & \text{otherwise. }
        \end{array}\right.$$
\end{lemma}

\medskip

\begin{remark}\label{symmetry_mm2}
If we use the original variables $x$ and $y$, where $y=1-z$, then the upper bound { function} for the difference $|C(x,y)-C(y,x)|$ for $x\le y$ is given by
$$\widehat{G}(x, y)=G(x,1-y) = \left\{\begin{array}{ll}
        \frac{x(y-x)}{y}; & \text{if } \max\{1-y,y^2\}\le x\le y, \\
        \frac{(1-y)(y-x)}{1-x}; & \text{if } x\le y \le\min\{2x-x^2,1-x\},  \\
        x(1-y); & \text{if } 1-y\le x\le y^2 \text{ or }2x-x^2\le y\le 1-x.
        \end{array}\right.$$
This follows directly from Lemma \ref{ocenaMM}.
By Remark \ref{symmetry_mm} the difference is symmetric also with respect to the line $x+y=1$. { So,} it follows that $\widehat{G}(x, y)=\widehat{G}(1-y,1-x)$ for all $x,y\in [0,1]$ with $x\le y$. The function $\widehat{G}$ will enable us to compute the maximal asymmetry function for the class of copulas $\mM^m$.

Observe also that functions $F(x,y)$ for the Marshall copulas and $\widehat{G}(x,y)$ for the maxmin copulas coincide in the region $\{(x,y);\ \max\{x,1-x\}\le y\le 1 \}$. Moreover,
$$\widehat{G}(x,y)=\left\{\begin{array}{ll}
       F(x,y); & \text{if } \frac12\le y\le 1,\  1-y\le x\le y,\\
        F(1-y,1-x); & \text{if } 0\le x\le\frac12,\ x\le y\le 1-x,
        \end{array}\right.$$
or $\widehat{G}(x,y)=\max\{F(x,y),F(1-y,1-x)\}$. This means that the maximal asymmetry functions for classes of copulas $\mM$ and $\mM^m$ coincide on the triangle $\{(x,y); \, x+y\le 1\}$. Furthermore, the maximal asymmetry function for $\mM^m$ is symmetric with respect to the line $x+y=1$, whereas the maximal asymmetry function for $\mM$ is not.
\qed
\end{remark}

\medskip

\begin{figure}[h]
\begin{center}
\includegraphics[width=7cm]{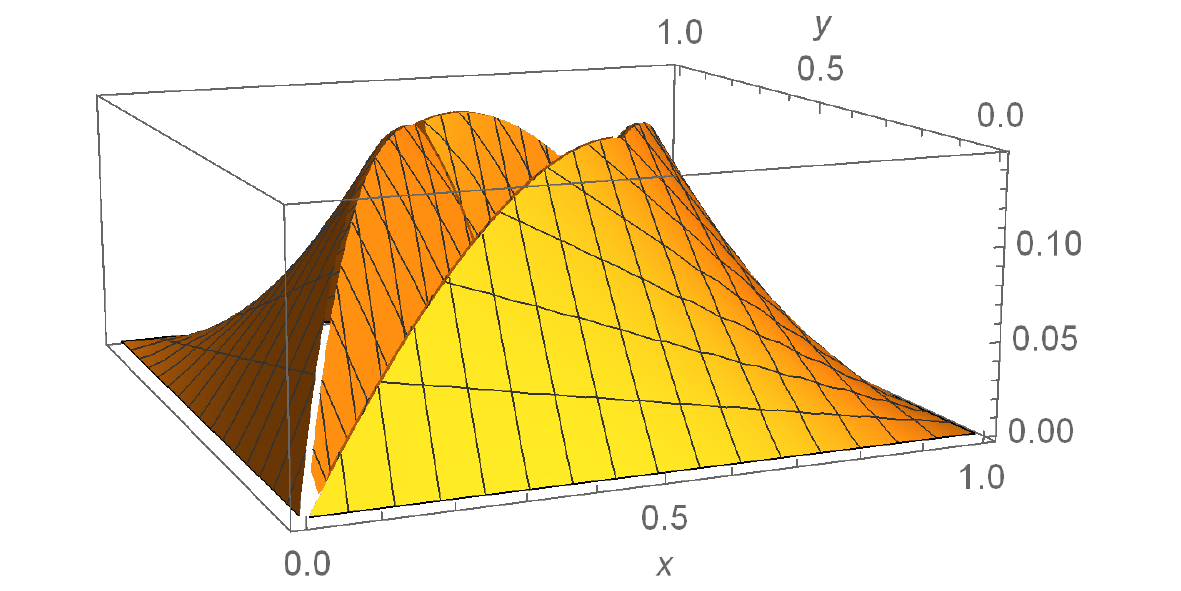} \quad \quad
\includegraphics[width=7cm]{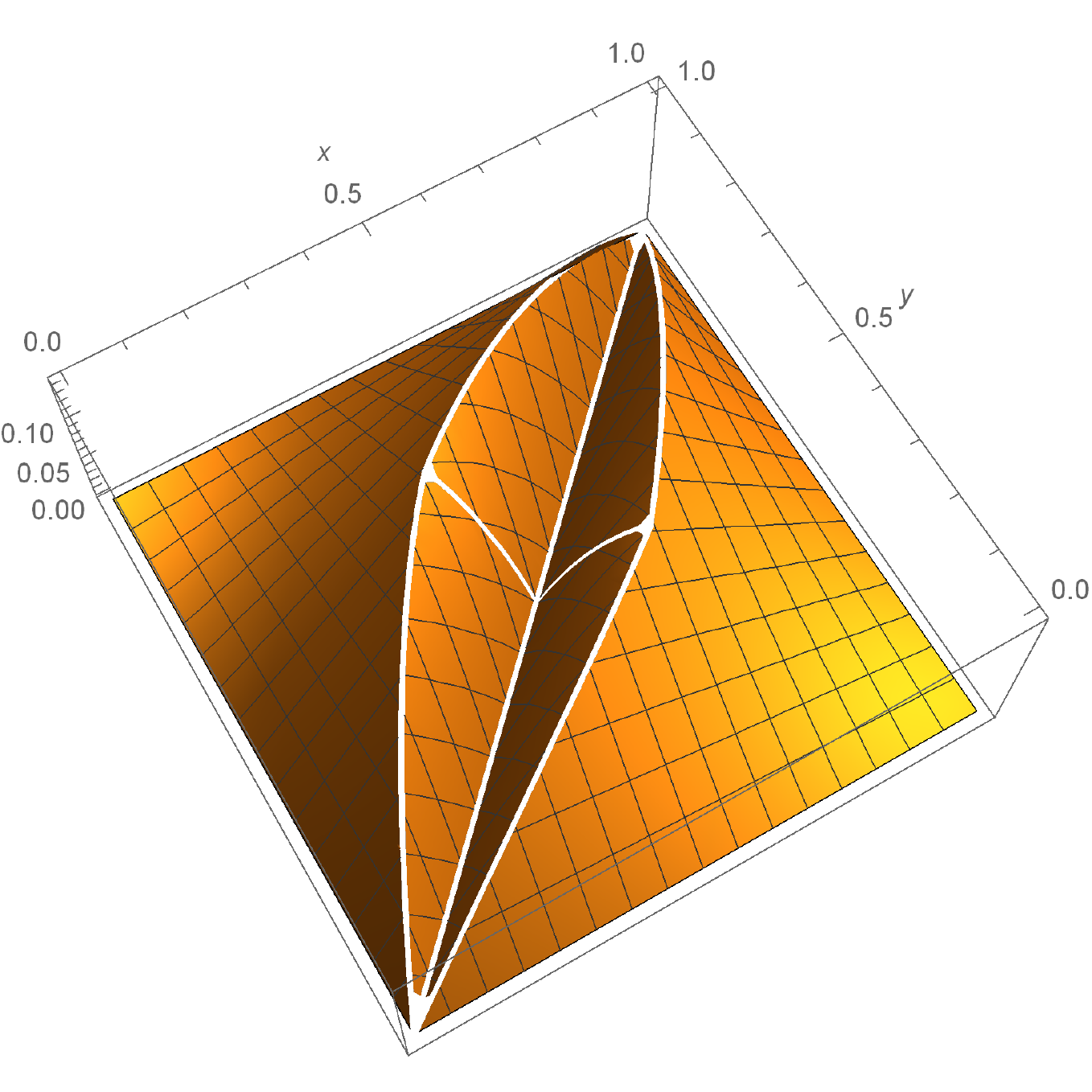}
\caption{Graph of $d^*_{\mM^m}$ for the set of maxmin copulas $\mM^m$.}\label{slika12}
\end{center}
\end{figure}

Figure \ref{slika12} shows two views of the upper bound { function}  $\widehat{G}$ for the expression $|C(x, y) - C(y, x)|$ for a maxmin copula $C$. We will see that this is in fact the maximal asymmetry function for the set of maxmin copulas. The point where the maximum of this function is attained is shown with a dot on the ridge of the graph.

\medskip

For $\mu\in (0,1)$, the maxmin copula defined by functions
$$ f_{\mu}(x) = \left\{\begin{array}{ll}
        \mu -x; & \text{if } 0<x\le \mu, \\
        0; & \text{if } x=0 \text{ or } \mu\le x \le 1,
        \end{array}\right.
\ \ \text{and} \ \
g(x) = \left\{\begin{array}{ll}
        0; & \text{if } x=0, \\
        1-x; & \text{if } 0 < x \le 1.
        \end{array}\right.$$
is equal to the Marshall copula $C_{\mu}$ of Example \ref{dosezenoM}. Next, let $D_{\mu}$, for $\mu\in (0,1)$, be maxmin copula defined by functions
$$ f(x) = \left\{\begin{array}{ll}
        0; & \text{if } x=0, \\
        1-x; & \text{if } 0 < x \le 1,
        \end{array}\right.
\ \ \text{and} \ \
g_{\mu}(x) = \left\{\begin{array}{ll}
        \mu -x; & \text{if } 0<x\le \mu, \\
        0; & \text{if } x=0 \text{ or } \mu\le x \le 1.
        \end{array}\right.
$$
So, we have
\begin{equation}\label{M_lambda}
    D_{\mu}(x, y) = \left\{\begin{array}{ll}
        x; & \text{if } 1-\mu+\mu x\le {y}\le 1,\\ %\text{ and } {y}\le \mu x+1-\mu, \\
        y-(1-\mu)(1-x); & \text{if } 1-\mu\le {y}\le 1-\mu + \mu x, \\
        xy; & \text{if } 0 \le y \le 1-\mu.
        \end{array}\right.
\end{equation}
Observe that $D_{\mu}(x,y)=x+y-1+C_{\mu}(1-y,1-x)$. { Consequently,}  $D_{\mu}$ is the survival copula of $C_{\mu}^t$ and we obtain
\begin{equation}\label{DeltaD_mu}
D_{\mu}(x,y)-D_{\mu}(y,x)=C_{\mu}(1-y,1-x)-C_{\mu}(1-x,1-y).
\end{equation}

\medskip

\begin{theorem}\label{d*_for_maxmin}
Let $\mM^m$ be the set of all maxmin copulas and $C\in \mM^m$. Then $\mu_\infty(C) \le \frac{4}{27}$ and the bound is attained, i.e., $\mu_{\infty}(\mM^m)=\frac{4}{27}$. The maximal asymmetry function for $\mM^m$ is equal to
$$d^*_{\mM^m}(x,y)= \left\{\begin{array}{ll}
        \frac{x}{y}(y-x); & \text{if } \max\{1-y,y^2\}\le x\le y, \\
        \frac{1-y}{1-x}(y-x); & \text{if } x\le y \le\min\{2x-x^2,1-x\},  \\
        x(1-y); & \text{if } y \ge \min\{ \sqrt{x}, 2x-x^2 \},\\
        y(1-x); & \text{if } x \ge \min\{ \sqrt{y}, 2y-y^2 \},\\
        \frac{y}{x}(x-y); & \text{if } \max\{1-x,x^2\}\le y\le x, \\
        \frac{1-x}{1-y}(x-y); & \text{if } y\le x \le\min\{2y-y^2,1-y\}.
        \end{array}\right.$$
\end{theorem}

\prf
By Example \ref{dosezenoM} copulas $C_{\mu}$ are extremal Marshall copulas, i.e., they are such that the upper bound $F(x,y)$ is attained  for each $(x,y)$ by at least one of them. Since $\widehat{G}(x,y)=F(x,y)$ for $1-y\le x\le y$ and since $C_{\mu}$ are also maxmin copulas, we have for at least one value of $\mu$ the equality $|C_{\mu}(x, y) - C_{\mu}(y,x)|=\widehat{G}(x,y)$ for each $(x,y)$ such that $1-y\le x\le y$.

Next, maxmin copulas $D_{\mu}$ are such that the upper bound $F(1-y,1-x)$ is attained  for each $(x,y)$ by at least one of them. (See Remark \ref{symmetry_mm2}, in particular equality (\ref{DeltaD_mu}).) Hence, we have for at least one value of $\mu$ the equality $|D_{\mu}(x, y) - D_{\mu}(y,x)|=\widehat{G}(x,y)$ for each $(x,y)$ such that $x\le y\le 1- x$.

Finally, recall from Example \ref{dosezenoM} that
$$\mu_\infty(C_{\frac{2}{3}}) = \frac{4}{27}.$$
For Marshall copulas we have $\frac{4}{27}=F(\frac49,\frac23)$. Hence, the maximal value of $\frac{4}{27}$ for $\widehat{G}$ is attained at $(\frac49,\frac23)$ and $(\frac59,\frac13)$. The corresponding copulas are $C_{\frac23}$ and $D_\frac23$.
\qed

\medskip

\begin{remark}
Since maxmin copulas are PQD we have $d^*_{\mM^m}(x,y)\le d^*_{\mP}(x,y)$ for all $(x,y)$. Observe that the two maximal asymmetry functions that are given in Theorems \ref{d*PN1} and \ref{d*_for_maxmin}, respectively, coincide for $(x,y)$ with $1-y\le x\le y^2$, or $2x-x^2\le y\le 1-x$, or $1-x \le y\le x^2$, or $2y-y^2\le x\le 1-y$.
\end{remark}

\medskip

\begin{observation} \label{zg_mejaMM_p}
Let $C$ be a maxmin copula and $p\in[1,\infty)$.
Then
\begin{equation}\label{mmp_ocena}
\mu_p(C) \le \left(2\int\!\!\int_{D_1} G(x,z)^p\,dx\,dz\right)^{1/p}
\end{equation}
where $D_1 = \{(x, z)\in [0, 1]^2, x+z \le 1\}$. In particular, $\mu_1(C) \le 0.0577$ and $\mu_2(C) \le 0.0707$.
If $C$ is a maxmin copula it follows from Lemma \ref{C-f-g} and Lemma \ref{ocenaMM} that
(\ref{mmp_ocena}) holds. The integral in (\ref{mmp_ocena}) cannot be computed analytically for general $p$. For $p=1, 2$ we get
$$ \mu_1(C) \le 2\int\!\!\int_{D_1} G(x,z)\,dx\,dz = $$
$$ = 2 \ln 2 - \textstyle \frac{2}{3}\ln(\sqrt{5}+3) + \frac{1}{36}(97 - 47 \sqrt{5}) \approx 0.0577,$$
$$ \mu_2(C) \le \left(2\int\!\!\int_{D_1} G(x,z)^2\,dx\,dz\right)^{1/2} = $$
$$ = \textstyle\frac{1}{12} \sqrt{2082 - \frac{3038}{3} \sqrt{5} + 864 \ln(\sqrt{5}-1)} \approx 0.0707.$$

We use (\ref{mmp_ocena}) to compute numerically an upper bound for asymmetry $\mu_p(\mM^m)$ for the set of maxmin copulas $\mM^m$. The value of $\mu_p(C_{\frac23})$ for the copula $C_{\frac23}$ is given in Example \ref{dosezenoM_p}.
Figure \ref{slika14} below shows the graphs of both { expressions} as functions of $p$ when $1 \le p \le 10$ (full curve and dashed curve), and also the asymptote of both, $y = \frac{4}{27}$ (dotted  line).
\end{observation}

\begin{figure}[h]
\begin{center}
\includegraphics[width=9cm]{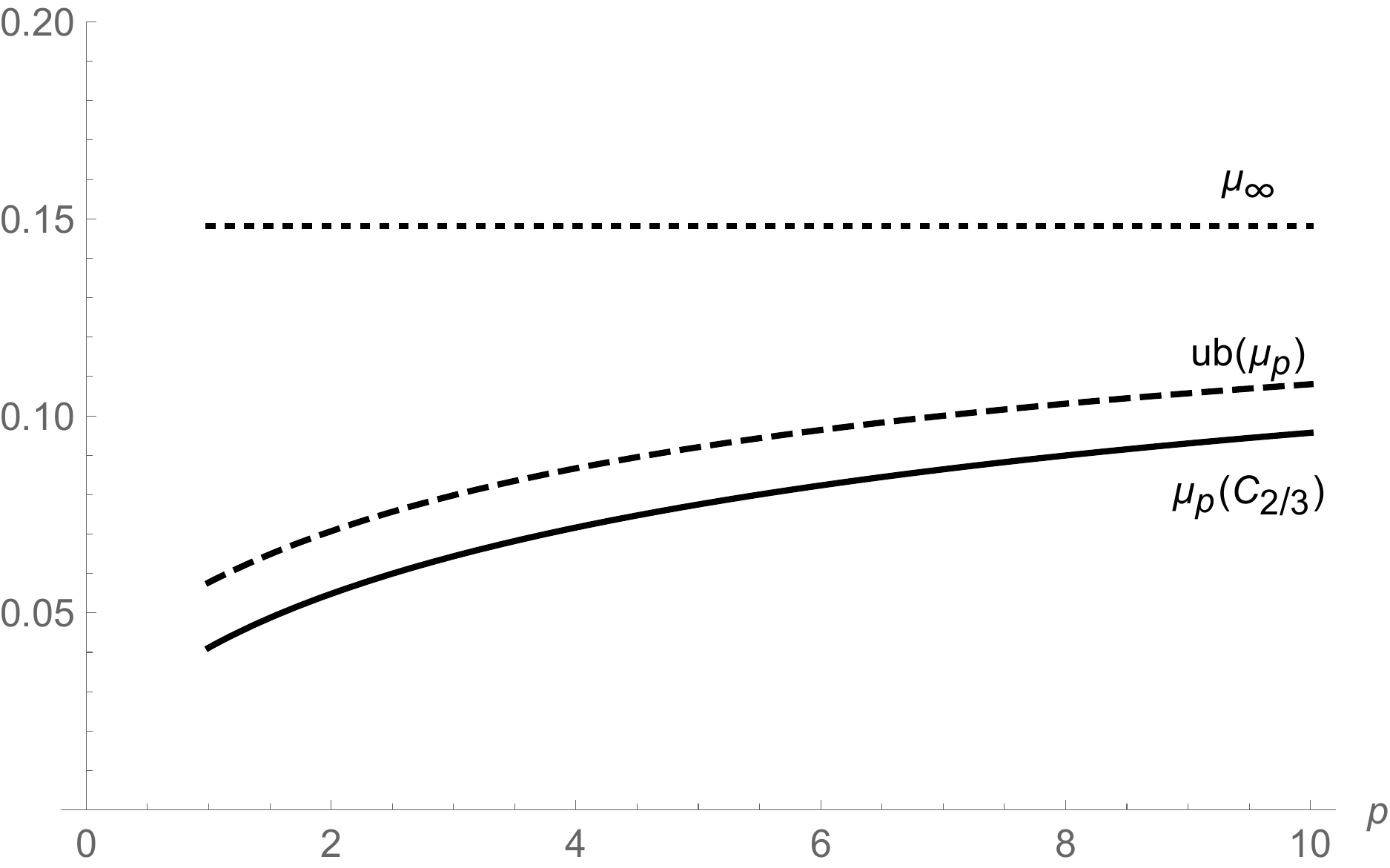}
\caption{Bounds for $\mu_p(\mM^m)$.}\label{slika14}
\end{center}
\end{figure}

\section{ Reflected maxmin copulas }

A reflected maxmin copula is {defined} by
$$C(x, y) = \max\{0, xy - f(x)g(y)\} ,$$
where $f, g: [0, 1] \to [0, 1]$ are functions in the set ${\mathcal F}$ introduced just before Lemma \ref{takojzaF}, see \cite{KO}.

Let $x, y \in [0, 1]$. Then
$$|C(x, y) - C(y, x)| = |\max\{0, xy - f(x)g(y)\} - \max\{0, xy - f(y)g(x)\}|.$$
Since the expression $|C(x, y) - C(y, x)|$ is symmetric with respect to the line $x=y$, we may assume that $y \ge x$.
It is also symmetric with respect to the interchange of functions $f$ and $g$, so we may assume that $f(x)g(y) \le f(y)g(x)$ for fixed $(x,y)\in [0,1]$.

\begin{lemma} \label{ocenaIMM}
Let $x, y \in [0, 1]$, $y \ge x$ and suppose $f(x)g(y) \le f(y)g(x)$. Then
$$|\max\{0, xy - f(x)g(y)\} - \max\{0, xy - f(y)g(x)\}| \le H(x,y)$$
where
$$H(x, y) = \left\{\begin{array}{ll}
        \frac{x(y-x)}{y}; & \text{if } y(1-y) \le x \le 1-y, \\
        \frac{(y-x)(1-y)}{y}; & \text{if } x \ge 1-y, \\
        xy; & \text{if } x \le y(1-y).
        \end{array}\right.$$
\end{lemma}

The proof of the lemma is rather lengthy and we include it in the Appendix.

\medskip

The {next} figure shows the upper bound for the expression $|C(x, y) - C(y, x)|$ for a reflected maxmin copula $C$. We will see that this is in fact the maximal asymmetry function for the set of reflected maxmin copulas.

\begin{figure}[h]
\begin{center}
\includegraphics[width=7cm]{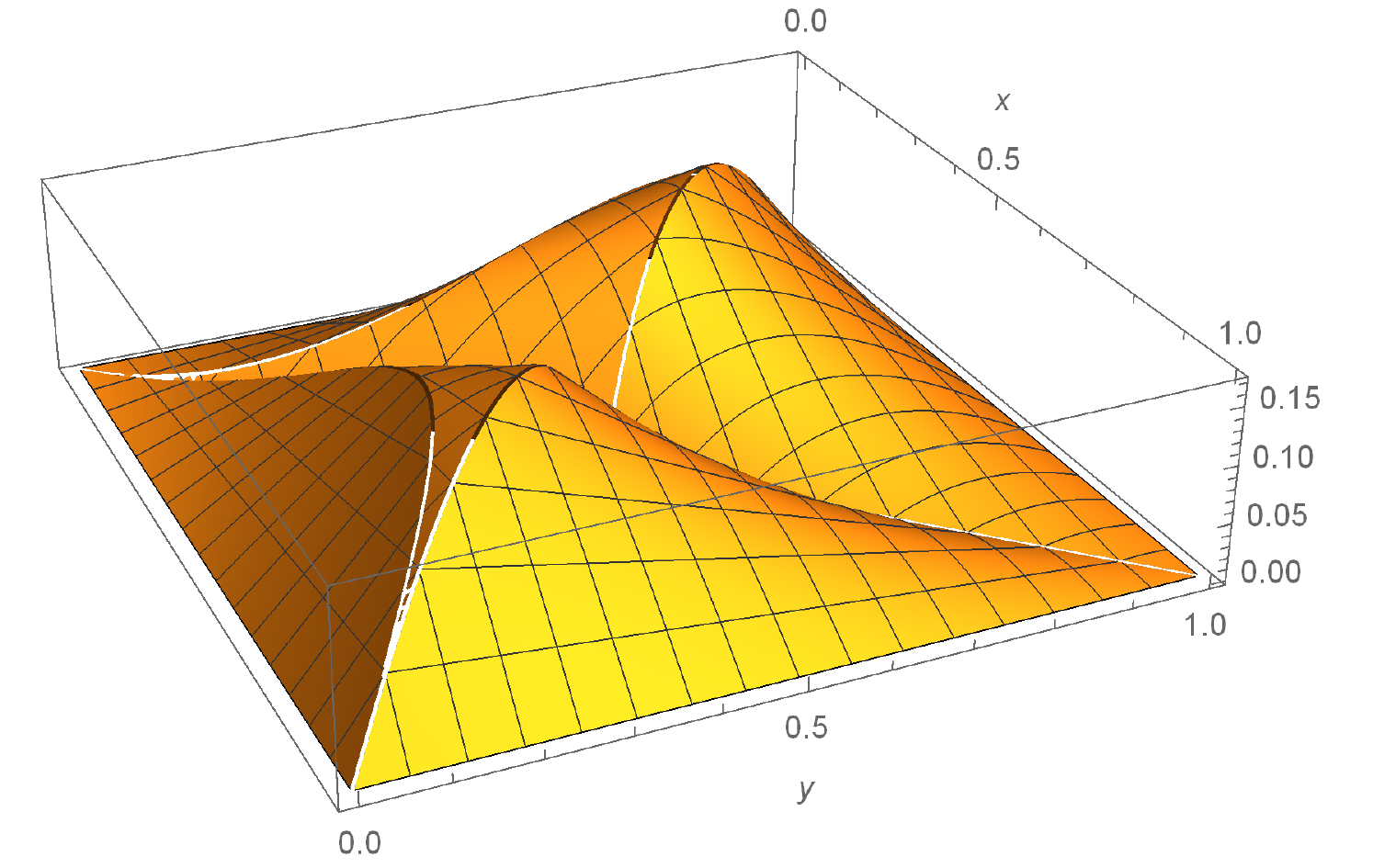} \quad \quad
\includegraphics[width=6cm]{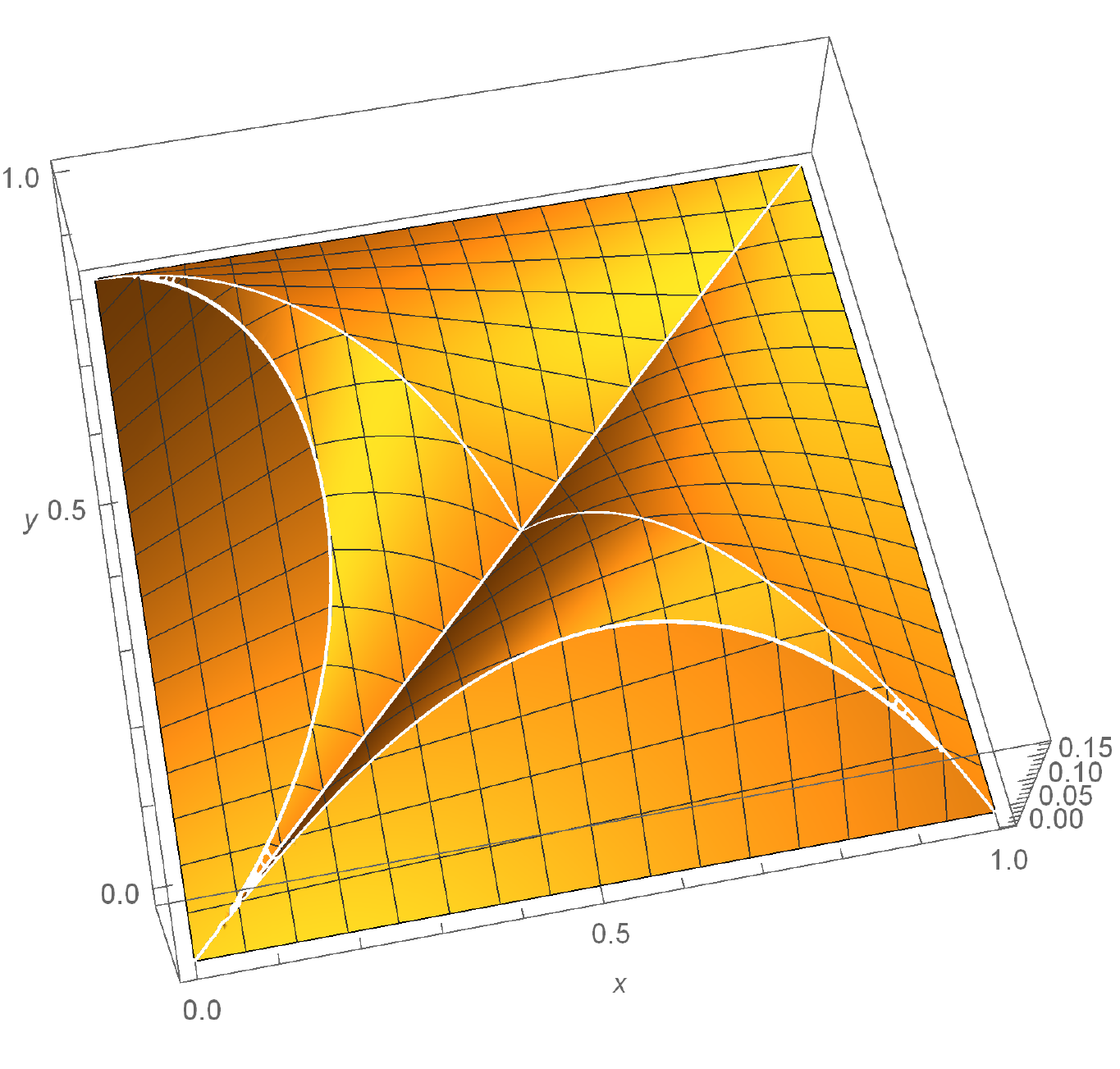}
\caption{Graph of $d^*_{\mM^{m\sigma}}$ of the set of reflected maxmin copulas $\mM^{m\sigma}$.}
\end{center}
\end{figure}

\begin{lemma} \label{zgornja_mejaIMM}
Let $x, y \in [0, 1]$ and $y \ge x$. Then
$$H(x, y) \le 3 - 2\sqrt{2}.$$
\end{lemma}

\prf
Function $H(x,y)$ is $0$ on the boundary of the triangle $x,y\in [0,1]$, $y\ge x$. Functions $\frac{x(y-x)}{y}$, $\frac{(y-x)(1-y)}{y}$ and $xy$ have no stationary points in this triangle, so the maximum has to be attained on the curves $x=y(1-y)$ or $x=1-y$. There we have $H(y(1-y),y)=y^2-y^3$ and $H(1-y,y)=\frac{(1-y)(2y-1)}{y}$. The first function has the maximum at $y=\frac23$. There we get $x=\frac29$ and $H(\frac29,\frac23)=\frac{4}{27}$. The second function has to be maximized on the interval $[\frac12,1]$. Its derivative is equal to $\frac{1-2y^2}{y^2}$, so the maximum is attained at $y=\frac{\sqrt{2}}{2}$. There we get $x=1-\frac{\sqrt{2}}{2}$ and $H(1-\frac{\sqrt{2}}{2},\frac{\sqrt{2}}{2}) = 3-2\sqrt{2}$, which is larger than $\frac{4}{27}$.
\qed

\medskip

Let $\lambda \in (0, 1], \mu \in [0, 1]$ be parameters and let $E_{\lambda,\mu}$ be a reflected maxmin copula defined by
functions $f_{\lambda}$ and
$g_{\mu}$, where
\begin{equation}\label{f_lam g_mu}
f_{\lambda}(x) = \left\{\begin{array}{ll}
        \frac{1-\lambda}{\lambda}x; & \text{if } x\le \lambda, \\
        1-x; & \text{if } x\ge \lambda,
        \end{array}\right.
				\ \ \text{and} \ \
g_{\mu}(x) = \left\{\begin{array}{ll}
        \mu - x; & \text{if } 0<x \le \mu, \\
        0; & \text{if }  x=0 \text{ or } \mu\le x \le 1.
        \end{array}\right.
\end{equation}
Then $f_{\lambda}, g_{\mu} \in {\mathcal F}$, and we have
\begin{equation}\label{E_lammu}
E_{\lambda,\mu}(x, y) = \left\{\begin{array}{ll}
        0; & \text{if } 0 \le y \le \min\{\mu(1-x),\mu(1-\lambda)\},  \\
        \frac{x}{\lambda}(y-\mu(1-\lambda)); & \text{if } 0\le x \le\lambda,\  \mu(1-\lambda)\le y\le \mu, \\
        xy; & \text{if } \mu\le y\le 1, \\
        \mu x+y-\mu; & \text{otherwise. }
        \end{array}\right.
\end{equation}

Observe in passing that copulas $E_{\lambda,\mu}$ are ordinal sums of copulas $\Pi$ and $W$, i.e., they are obtained via patchwork construction of copulas. This was already mentioned for the family $C_{\mu}$ in Section \ref{Marshall_copulas} and holds also for the family $D_{\mu}$, since $D_{\mu}$ is a survival copula of $C_{\mu}$. We refer to \cite[pp. 63-64]{Nelsen} for the notion of an ordinal sum of copulas and to \cite{DS-PS} for the general patchwork construction for bivariate copulas.

\medskip

\begin{theorem}\label{d*_RMM}
Let $\mM^{m\sigma}$ be the set of all reflected maxmin copulas and $C\in \mM^{m\sigma}$. Then $\mu_\infty(C) \le 3-2\sqrt{2}$ and the bound is attained, i.e., $\mu_{\infty}(\mM^{m\sigma})=3-2\sqrt{2}$. The maximal asymmetry function for $\mM^{m\sigma}$ is equal to
$$d^*_{\mM^{m\sigma}}(x,y)= \left\{\begin{array}{ll}
        \frac{x}{y}(y-x); & \text{if } y(1-y)\le x\le \min\{y,1 - y\}, \\
        \frac{1-y}{y}(y-x); & \text{if } 1-y\le x\le y,  \\
        xy; & \text{if } x\le y(1-y) \text{ or }  y\le x(1-x),\\
        \frac{y}{x}(x-y); & \text{if } x(1-x)\le y\le \min\{x,1 - x\}, \\
        \frac{1-x}{x}(x-y); & \text{if } 1-x\le y\le x.
        \end{array}\right.$$
\end{theorem}

\prf
We may assume that $x \le y$. We consider first a point $(x,y)$ in the triangle $$R_1=\textstyle\{(x,y);\ \frac12\le y\le 1,\ 1-y\le x\le y\}$$ and
choose $\lambda = y$ and $\mu = 1$. Observe that $E_{y,1}(x,y) = \frac{x}{y}(2y-1)$ and $E_{y,1}(y,x) = x+y-1$, so $E_{y,1}(x,y)-E_{y,1}(y,x) = \frac{1-y}{y}(y-x)$.
It follows that the upper bound $H(x,y)$ is attained by the family $E_{\lambda,1}$ on this triangle. Note also, that we have $H\left(1-\frac{\sqrt{2}}{2},\frac{\sqrt{2}}{2}\right) = 3-2\sqrt{2}$ and that the point $\left(1-\frac{\sqrt{2}}{2},\frac{\sqrt{2}}{2}\right)$ is in $R_1$.
Therefore,
$$\mu_\infty\left(E_{\sqrt{2}/2,1}\right) = 3 - 2\sqrt{2}$$
and the upper bound of Lemma \ref{zgornja_mejaIMM} is attained. Next, we consider a point $(x,y)$ in the triangle
$$R_2=\textstyle\{(x,y);\ 0\le x\le \frac12,\ x\le y\le 1-x\}$$ and
choose $\lambda = y$ and $\mu = \frac{x}{1-y}$, which lies in the interval $[0,1]$. Observe that
$$E_{y,\frac{x}{1-y}}(x,y) = \left\{\begin{array}{ll}
        \frac{x}{y}(y-x); & \text{if } x \ge y(1-y),\ x\le y\le 1-x,\\
        xy; & \text{if } x \le y(1-y) ,
        \end{array}\right.$$
and $E_{y,\frac{x}{1-y}}(y,x) =0.$
It follows that the upper bound $H(x,y)$ on this triangle is attained by the family $E_{\lambda,\mu}$.
\qed

\medskip

\begin{figure}
\begin{center}
\includegraphics[width=4cm]{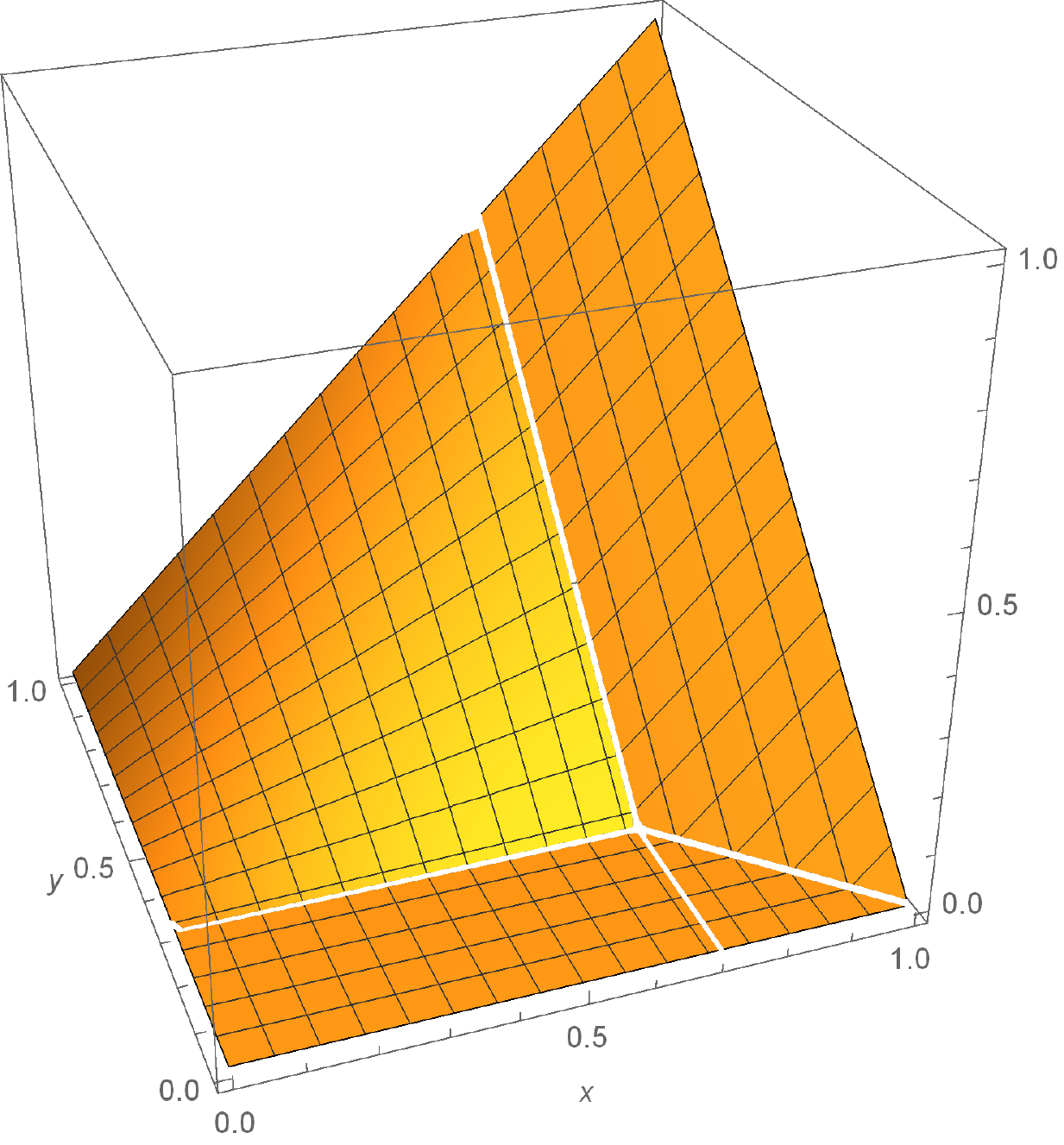} \ \
\includegraphics[width=42mm]{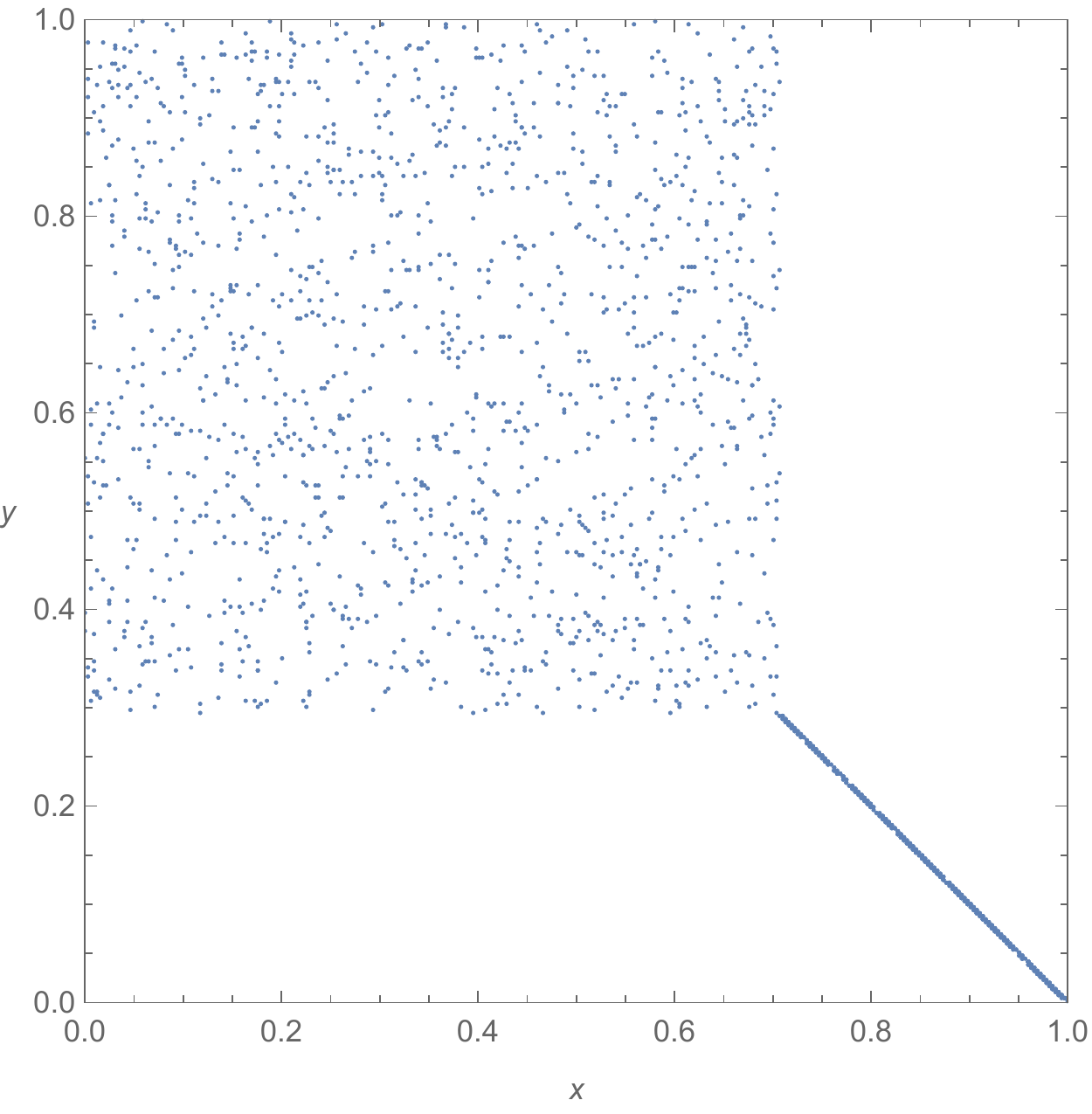} \ \
\includegraphics[width=6cm]{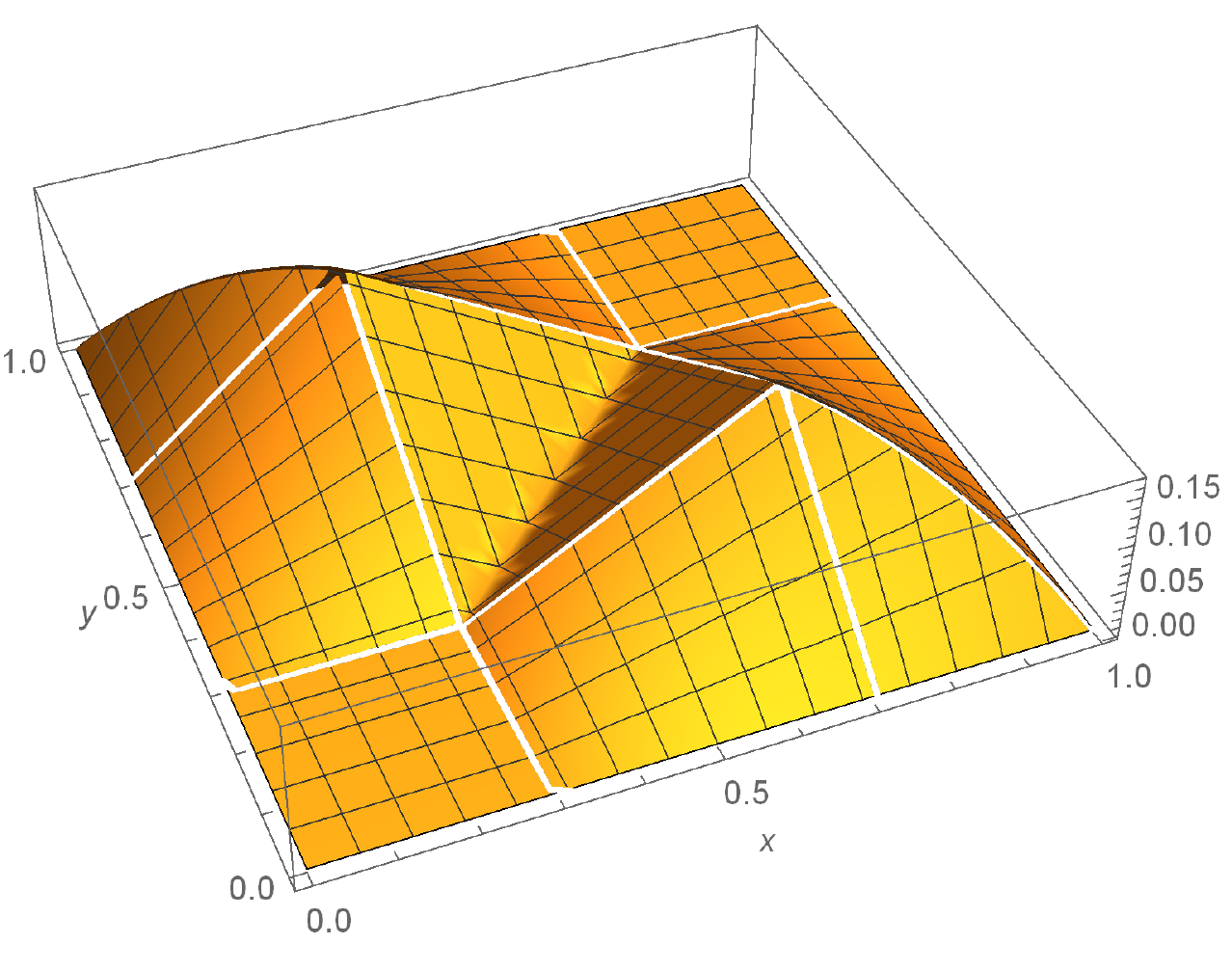}
\caption{Figure shows the graph of copula $E_{\sqrt{2}/2,1}$, its scatterplot of 2000 points, and the graph of the expression $|E_{\sqrt{2}/2,1}(x, y) - E_{\sqrt{2}/2,1}(y, x)|$, respectively.}\label{slika19}
\end{center}
\end{figure}

\medskip

\begin{remark}
Since reflected maxmin copulas are NQD we have $d^*_{\mM^{m\sigma}}(x,y)\le d^*_{\mN}(x,y)$ for all $(x,y)$. Observe that the two maximal asymmetry functions that are given in Theorems \ref{d*PN2} and \ref{d*_RMM}, respectively, coincide for $(x,y)$ with $x\le y(1-y)$ or $y\le x(1-x)$.
\end{remark}

\medskip

If $C$ is a reflected maxmin copula it follows from Lemma \ref{ocenaIMM} that
$$ \mu_p(C) \le \left(2\int\!\!\int_D H(x,y)^p\,dx\,dy\right)^{1/p} $$
where $D$ is triangle $D = \{(x, y)\in [0, 1]^2, y \ge x\}$.
This integral cannot be computed analytically for general $p$. For $p=1, 2$ we get
$$ \mu_1(C) \le 2\int\!\!\int_{D} H(x,y)\,dx\,dy = \textstyle \frac{1}{3}\ln 2 - \frac{31}{180} \approx 0.0588,$$
$$ \mu_2(C) \le \left(2\int\!\!\int_{D} H(x,y)^2\,dx\,dy\right)^{1/2} = \textstyle\frac{1}{6} \sqrt{84 \ln 2 - \frac{2437}{42}} \approx 0.0746.$$

\begin{example} \label{dosezeno_pIMM}
Let $E_{\sqrt{2}/2,1}$ be  the maxmin copula defined above and $p \ge 1$.
Then $$\mu_p(E_{\sqrt{2}/2,1}) = \left(\frac{2B(\sqrt{2}-1; p+2, p+1)}{\sqrt{2}^p(p+1)} + \frac{2(3-2\sqrt{2})^{p+1}((\sqrt{2}+1)p+2\sqrt{2}+1)}{(p+1)^2(p+2)}\right)^{1/p}.$$
In particular, $$\mu_1(E_{\sqrt{2}/2,1}) = \textstyle\frac{29}{24}\sqrt{2} -  \frac{5}{3} \approx 0.0422$$ and
$$\mu_2(E_{\sqrt{2}/2,1}) = \textstyle\sqrt{\frac{71}{45}\sqrt{2} -  \frac{401}{180}}\approx 0.0595.$$

We integrated the function $|E_{\sqrt{2}/2,1}(x,y)-E_{\sqrt{2}/2,1}(y,x)|$ over the triangle $T$ using the Mathematica software \cite{Mathematica} and obtained the value given above.

Thus, we have found an upper bound function of $p$ for asymmetry $\mu_p(C)$  in terms of the incomplete Beta function for a maxmin copula $C$  by taking $C=E_{\sqrt{2}/2,1}$.
\end{example}

\medskip

\begin{figure}[h]
\begin{center}
\includegraphics[width=9cm]{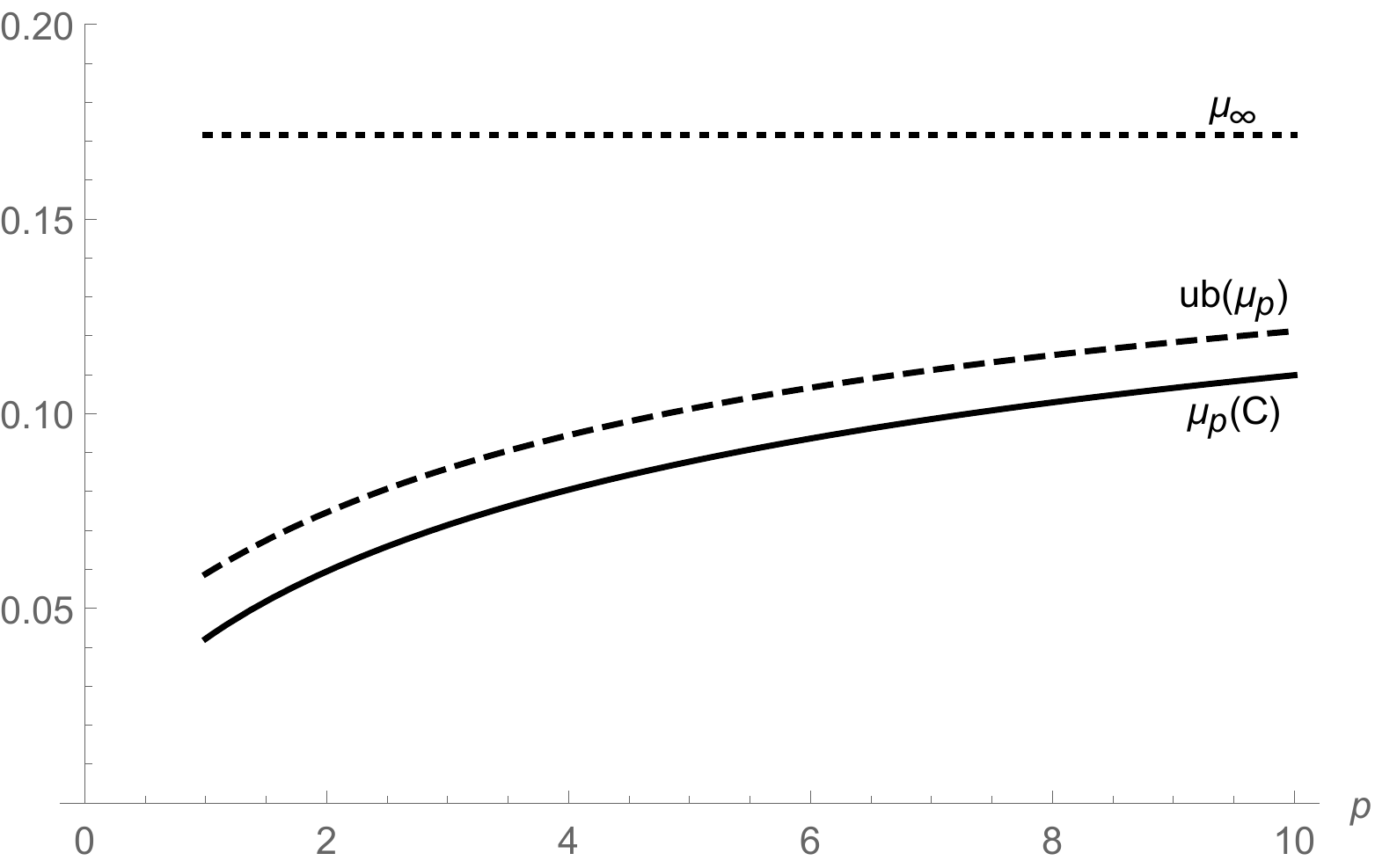}
\caption{Bounds for $\mu_p(\mM^{m\sigma})$.}\label{slika22}
\end{center}
\end{figure}

Figure \ref{slika22} presents  graphs of  $\mu_p(C)$ and the upper bound function as functions of $p$ when $1 \le p \le 10$ (full curve and dashed curve respectively), and also the asymptote of both functions, $y = 3-2\sqrt{2}$ (dotted line).

\section{Stochastic interpretation of the extremal families}

In this section we give stochastic interpretation of those families of copulas presented in the previous three sections for which the extremal values of the maximal asymmetry functions are attained. That is, the extremal families of copulas were constructed in such a way that the maximal possible difference $|C(a,b)-C(b,a)|$ that is allowed at a given point $(a,b)$ for the particular family of copulas was attained. It turned out that a copula constructed in such a way is extremal not only at the point $(a,b)$ but on a segment containing that point.

In the construction of shock models that correspond to the extremal copulas for each of the three families of copulas in shock models considered, we proceed as follows: First, we recall the relations that hold among the idiosyncratic shocks $X$ and $Y$ and the systemic shock $Z$ on one side and the resulting random vector $(U,V)$ on the other. These relations are specific for each of the three families. The generators of an extremal copula then govern how the three shocks are applied in a sequence. For instance, if Marshall copula $C_{\mu}$ governs the dependence structure of a shock model then we have that $P[Y\le Z]=1$ and $P[X\le Z]=\mu$, i.e., the idiosyncratic shock on the second component almost certainly precedes the systemic shock, while the probability that the idiosyncratic shock on the first component precedes the systemic shock is equal to the value of parameter $\mu$. Similar interpretations are obtained also for families $C_{\mu}$ and $D_{\mu}$ for maxmin copulas and family $E_{\lambda,\mu}$ of RMM copulas. Let us stress here that the interpretation of copulas $C_{\mu}$ as the extremal copulas of the family of maxmin copulas is quite different from the interpretation given above when copulas $C_{\mu}$ are considered as Marshall copulas. These come from very different relations among the shocks $X$, $Y$ and $Z$ in each of the two families of copulas.\\

We begin with family $C_{\mu}$, $\mu\in (0,1)$, of Marshall copulas presented in Example \ref{dosezenoM}. { We know from Theorem \ref{d* for Marshall} that these}   copulas give the maximal values of asymmetry for Marshall copulas.
{ Let us recall some general properties that Marshall copulas possess (see \cite{M}). We denote by $F_X$ the distribution function of a random variable $X$ and by $\supp(F_X)$ its support, i.e.\ the complement of the union of all open intervals $(a,b)\subset\mathbb{R}$ with $F_X(b)= F_X(a)$.}

Suppose that random vector $(U,V)$ has the joint distribution function given by
$H(u,v)=C(F_U(u),F_V(v)),$ where $C(u,v)=\min\{ug(v),f(u)v\}$ is a Marshall copula.  By \cite[Proposition  3.2]{M}, there are independent random variables $X$, $Y$ and $Z$ such that
$$U=\max\{X,Z\}\ \text{and}\ V=\max\{Y,Z\},$$ and $H(u,v)=F_X(u)F_Y(v)F_Z(\min\{u,v\})$. Furthermore, we have that
\begin{equation}\label{FFrelations}
F_U(x)=F_X(x)F_Z(x),\ F_V(y)=F_Y(y)F_Z(y),\text{and}\ f(F_U(x))\, F_V(x)=g(F_V(x))\, F_U(x).
\end{equation}
If $F_Z(x)\neq 0$, resp. $F_Z(y)\neq 0$, it follows that
\begin{equation}\label{FFrelations2}
f(F_U(x))=F_X(x)=\frac{F_U(x)}{F_Z(x)}\ \text{resp.}\  g(F_V(y))=F_Y(y)=\frac{F_V(y)}{F_Z(y)}.
\end{equation}

Consider now copula $C_{\mu}$ of (\ref{C_mu}). Since $g(t)=1$ for $t\in (0,1]$, the second relations of (\ref{FFrelations}) and (\ref{FFrelations2}) imply that $y\le z$ for any $y\in\supp(Y)$ and $z\in\supp(Z)$. Hence, we have that $\supp(Y)\subset (-\infty,z_1]$ and $\supp(Z)\subset [z_1,\infty)$ for some $z_1\in\mathbb{R}$. Since $f_{\mu}(t)=\max\{\mu,t\}$ for $t\in (0,1]$, the first relations of (\ref{FFrelations}) and (\ref{FFrelations2}) imply that $\supp(X)\subset (-\infty,z_1]\cup [z_2,\infty)$ and $\supp(Z)\subset [z_1,z_2]$, where $ z_1\le z_2$ and $F_X(z_1)=P[X\le z_1]=\mu$ and $P[X\ge z_2]=1-\mu$. Hence, $P[Y\le Z]=1$ and $P[X\le Z]=\mu$, i.e., it is almost sure that the shock $Y$ occurs before shock $Z$ and the probability that $X$ occurs before $Z$ is equal to the value of parameter $\mu$. The distribution functions of $U$ and $V$ are then equal to
$$F_U(x)= \left\{\begin{array}{ll}
%        0, & \text{if }  x\le z_1,\\
        \mu F_Z(x); & \text{if } x\le z_2,\\
        F_X(x); & \text{if }  x\ge z_2.
        \end{array}\right.\ \text{and}\ F_V(y)=F_Z(y).$$
The maximal possible asymmetry $\mu_{\infty}(\mM)=\frac{4}{27}$ is attained for $\mu=\frac23$.

\begin{figure}[h]
\begin{center}
\includegraphics[width=9cm]{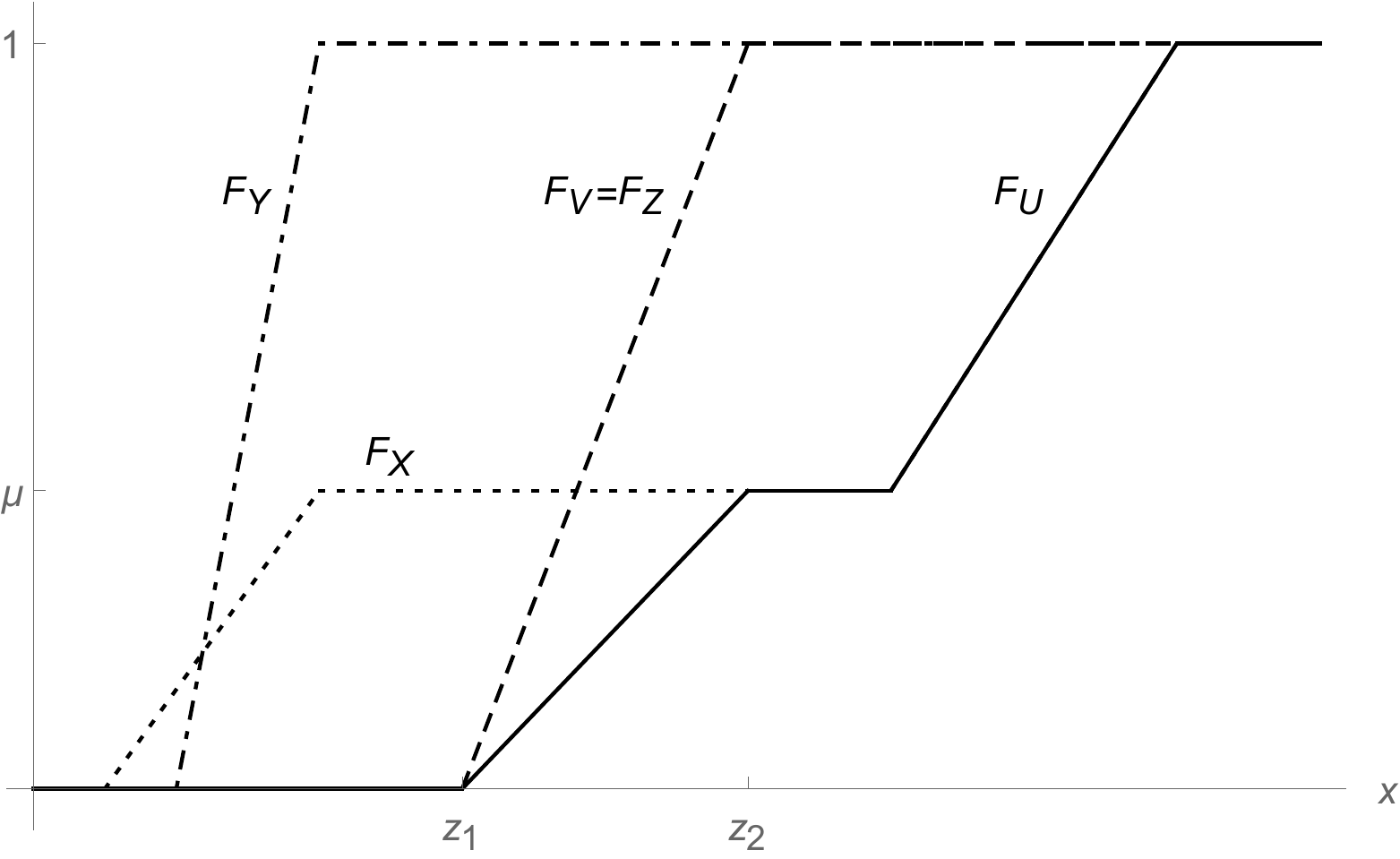}
\caption{Stochastic interpretation of $C_{\mu}$ as a Marshall copula.}\label{slika11}
\end{center}
\end{figure}

Graphs of possible distribution functions of shocks $X$, $Y$ and $Z$ are illustrated in Figure \ref{slika11}. {It clearly shows the interval $[z_1,z_2]$ where shock $Z$ is acting and $F_U=\mu F_Z$, while $X$ is acting partly before and partly after this interval. In the same figure we point out that $F_V=F_Z$ using a dashed line.}

\medskip

We continue this section with some comments on stochastic interpretation of the families of maxmin copulas $C_{\mu}$ and $D_{\mu}$, $\mu\in [0,1]$, { defined by}  (\ref{C_mu}) and (\ref{M_lambda}), respectively. These copulas give the maximal values of asymmetry for maxmin copulas described in Theorem \ref{d*_for_maxmin}.

Let us recall some general properties { of}  maxmin copulas  \cite{OR}.
Suppose that random vector $(U,V)$ has the joint distribution function given by $H(u,v)=C(F_U(u),F_V(v))$, where $C(u,v)$ is a general maxmin copula given in (\ref{maxmin}). By \cite[Theorem  10]{OR}, there are independent random variables $X$, $Y$ and $Z$ such that
$$U=\max\{X,Z\}\ \text{and}\ V=\min\{Y,Z\},$$
and by \cite[Lemma  7]{OR}, we have
\begin{equation}\label{FFrelationsMM0}
F_U(x)=F_X(x)F_Z(x),\ F_V(y)=F_Y(y)+F_Z(y)-F_Y(y)F_Z(y)
\end{equation}
{and so}
\begin{equation}\label{FFrelationsMM}
1-F_V(y)=(1-F_Y(y))(1-F_Z(y)).
\end{equation}
Assume that $F_U(x)>0$ and $F_V(x)<1$. Then the generators $\phi$ and $\psi$ and the distribution functions $F_U$ and $F_V$ are related by
\begin{equation*}\label{FFrelations2MM}
\phi(F_U(x))\left(F_V(x)-\psi(F_V(x))\right)=F_U(x)\left(1-\psi(F_V(x))\right)
\end{equation*}
(See \cite[Eq. (7)]{OR}.) Under the same assumptions it also holds that
\begin{equation}\label{FFrelations3MM}
\phi(F_U(x))=F_X(x)\ \text{and}\  \psi(F_V(y))=F_Y(y).
\end{equation}
For these examples we assume that  shocks $X$, $Y$ and $Z$ all have continuous distribution functions.\\

Consider now $C_{\mu}$ of (\ref{C_mu}) as a maxmin copula. Then
$$\phi_{\mu}(t)=\left\{\begin{array}{ll}
        0; & \text{if } t=0, \\
        \mu; & \text{if } 0\le t\le \mu, \\
        t; & \text{if } t\ge \mu,
        \end{array}\right.  \ \text{and}\ \psi(t)=\left\{\begin{array}{ll}
        0; & \text{if } t<1, \\
        1; & \text{if } t=1.
        \end{array}\right.$$
Note that the range of $\phi_{\mu}$ contains $0$ and the interval $[\mu,1]$. Then, if $x\in\mathbb{R}$ is such that $0<F_X(x)<\mu$ the first relations of (\ref{FFrelationsMM0}) and (\ref{FFrelations3MM}) imply that $F_Z(x)=0$. On the other hand, if $x\in\mathbb{R}$ is such that $\mu<F_X(x)\le 1$ then the same relations imply that $F_Z(x)=1$. Therefore it follows that $\supp(X)\subset (-\infty,z_1]\cup [z_2,\infty)$, $\supp(Z)\subset [z_1,z_2]$ for some $z_1$ and $z_2$ with $z_1\le z_2$, and $P[X\le z_1]=\mu$ and $P[X\ge z_2]=1-\mu$. Furthermore, we have that
$$F_U(x)=\left\{\begin{array}{ll}
        \mu F_Z(x); & \text{if } x\le z_2, \\
        F_X(x); & \text{if } x\ge z_2.
        \end{array}\right.$$
Since the range of $\psi$ contains only $0$ and $1$  the second relation of (\ref{FFrelations3MM}) implies that if $F_V(x)<1$ then $F_Y(x)=0$. Next relation (\ref{FFrelationsMM}) implies that if $F_V(x)<1$ we have that $F_V(x)=F_Z(x)$. If $x_1$ is minimal such that $F_V(x_1)=1$ then $F_Y(x)=0$ for all $x<x_0$. Since $Y$ is a continuous random variable also $F_Y(x_0)=0$. Then (\ref{FFrelationsMM}) implies that $F_Z(x_0)=1$. Hence $F_V(x)=F_Z(x)$ for all $x$.

\begin{figure}[h]
\begin{center}
\includegraphics[width=9cm]{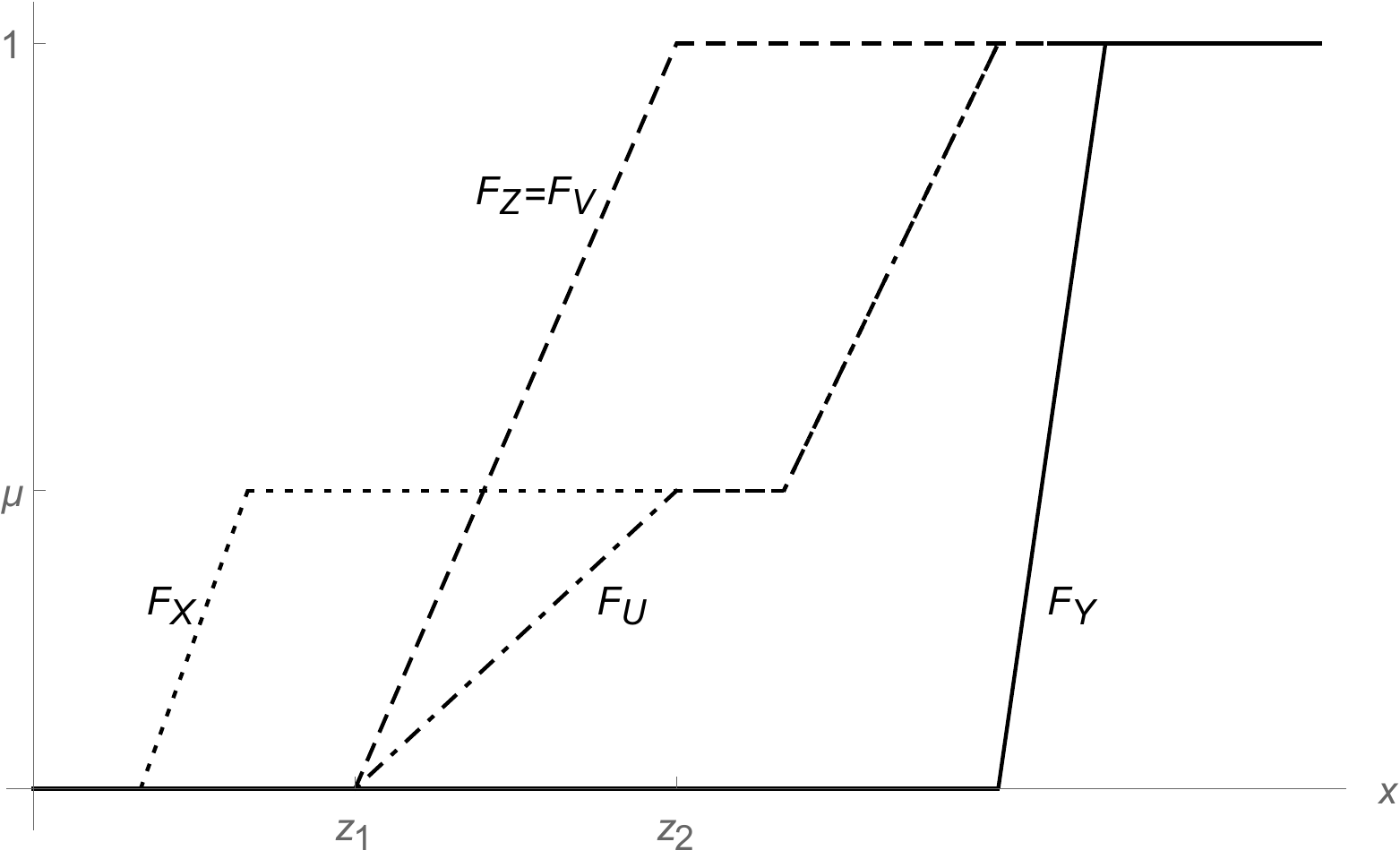}
\caption{Stochastic interpretation of $C_{\mu}$ as a maxmin copula.}\label{slika15}
\end{center}
\end{figure}

Graphs of possible distribution functions of shocks $X$, $Y$ and $Z$ are illustrated in Figure \ref{slika15}.
It clearly shows the interval $[z_1,z_2]$ where shock $Z$ is acting and $F_U=\mu F_Z$, while $X$ is acting partly before and partly after this interval. In the same figure we point out that $F_V=F_Z$ using a dashed line. In difference with Figure \ref{slika11} shock $Y$ is here acting after shock $Z$.

\begin{figure}[h]
\begin{center}
\includegraphics[width=9cm]{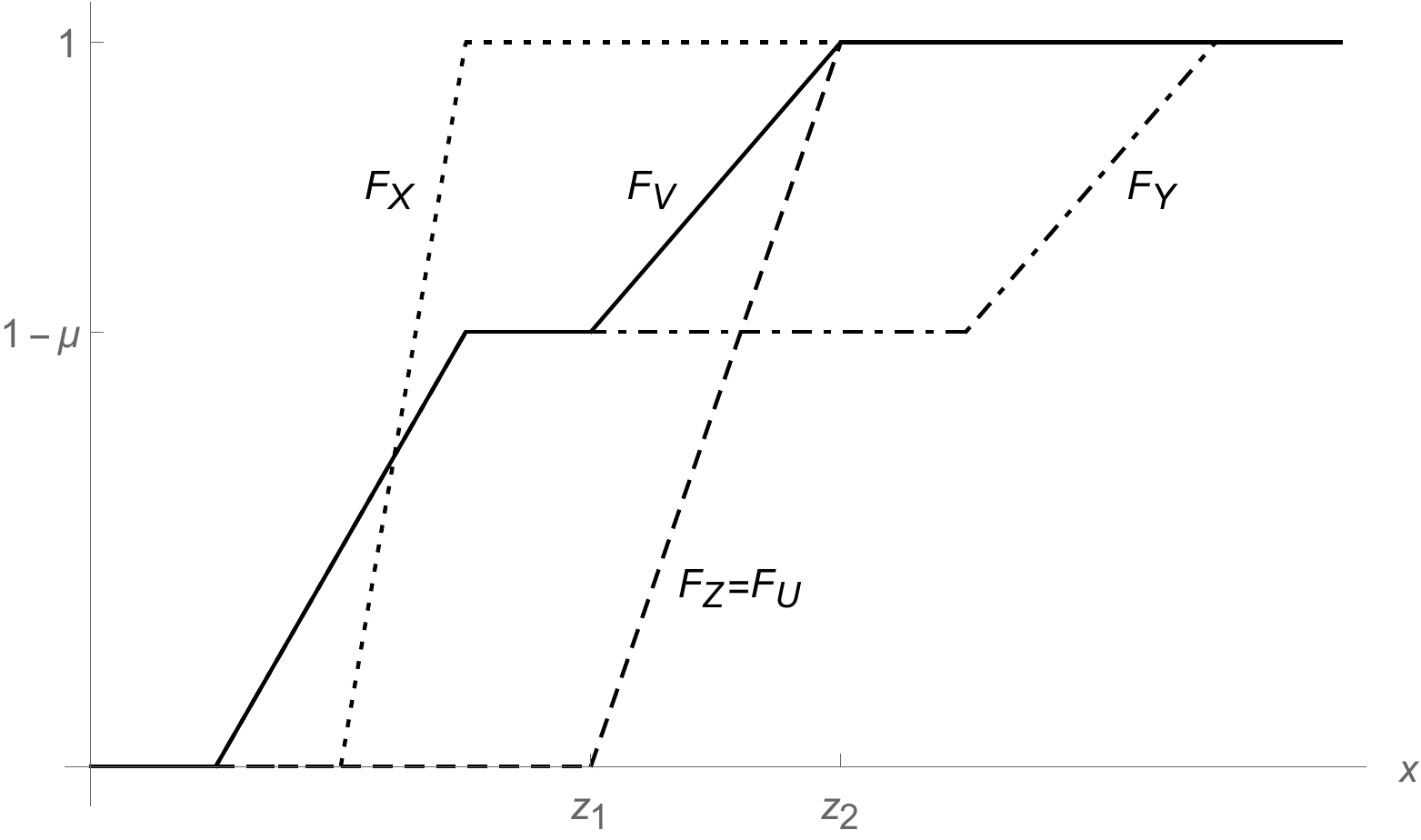}
\caption{Stochastic interpretation of maxmin copula $D_{\mu}$.}\label{slika16}
\end{center}
\end{figure}

Consider next the maxmin copula $D_{\mu}$ of (\ref{M_lambda}). Its generators are
$$\phi(t)=\left\{\begin{array}{ll}
        0; & \text{if } t=0, \\
        1; & \text{if } t>0,
        \end{array}\right. \ \text{and}\ \psi_{\mu}(t)=\left\{\begin{array}{ll}
        t; & \text{if } 0\le t\le 1-\mu, \\
        1-\mu; & \text{if }  1-\mu\le t <1,\\
	 1; & \text{if } t=1.
        \end{array}\right.$$
Since $\phi(t)=1$ for $t\in (0,1]$ the first relations of (\ref{FFrelationsMM0}), %relation (\ref{FFrelationsMM}), and relations
and  (\ref{FFrelations3MM}) imply that $x\le z$ for any $x\in\supp(X)$ and $z\in\supp(Z)$. Hence, we have that $\supp(X)\subset (-\infty,x_2]$ and $\supp(Z)\subset [z_1,\infty)$ for some $x_2$ and $z_1$ with $x_2\le z_1$ and $P[X\le Z]=1$.
Therefore, we have that $$F_U(x)=F_Z(x).$$
Relation (\ref{FFrelationsMM}) and the second relation (\ref{FFrelations3MM}) imply
$$1-\psi_{\mu}\left(1-\left(1-F_Y(y)\right)\left(1-F_Z(y)\right)\right)=1-F_Y(y).$$
Since $\psi_{\mu}(t)=\min\{1-\mu,t\}$ for $t\in (0,1]$ it follows that $1-\psi_{\mu}(t)=\max\{\mu,1-t\}$. This implies that $\supp(Y)\subset (-\infty,y_2]\cup [y_3,\infty)$ and $\supp(Z)\subset [z_1,z_2]$, where $y_2\le z_1\le z_2\le y_3$ and $F_Y(y_2)=P[Y\le y_2]=1-\mu$ and $P[Y\ge y_3]=\mu$. Relation (\ref{FFrelationsMM}) then implies that
$$F_V(y)=\left\{\begin{array}{ll}
        F_Y(y); & \text{if } y<z_1, \\
        1-\mu+\mu F_z(y); & \text{if } z_1\le y\le z_2,\\
        1; & \text{if } y\ge z_2.
        \end{array}\right.$$
Graphs of possible distribution functions of shocks $X$, $Y$ and $Z$ are illustrated in Figure \ref{slika16}. It clearly shows the interval $[z_1,z_2]$ where shock $Z$ is acting and $F_U=\mu F_Z$, while $X$ is acting before this interval, meanwhile $Y$ is acting partly before and partly after this interval. In the same figure we trace out that graph of $F_V$ using a full line.

\medskip

To conclude let us comment on stochastic interpretation of the family of RMM copulas $E_{\lambda,\mu}$, for $\lambda$, $\mu\in (0,1)$, defined by (\ref{E_lammu}). These copulas give the maximal values of asymmetry for reflected maxmin copulas described in Theorem \ref{d*_RMM}.
 We first recall some general properties  of RMM copulas from \cite{KO}.

Suppose that random vector $(U,V)$ has the joint distribution function given by $H(u,v)=C(F_U(u),F_V(v))$, where $C(u,v)=\max\{0,uv-f(u)g(v)\}$ is a reflected maxmin copula.  By \cite[Theorem  19]{KO}, there are independent random variables $X$, $Y$ and $Z$ such that $U=\max\{X,Z\}$, $W=\min\{Y,Z\}$ and the distribution function of e.g. $V=-W$ is given by the survival function of $W$: we have $F_V(v)=1-F_W(-v)$ in the case of continuous distribution functions. By the relationship of maxmin copulas and reflected maxmin copulas given in \cite[Section  2]{KO}, we then have that
\begin{equation}\label{FFrelationsRMM}
F_U(x)=F_X(x)F_Z(x),\ \text{and}\ 1-F_V(-y)=F_Y(y)+F_Z(y)-F_Y(y)F_Z(y)\ \text{or}\ F_V(-y)=(1-F_Y(y))(1-F_Z(y)).
\end{equation}
For $\phi(t)=t+f(t)$ and $\psi(t)=t-g(1-t)$ it also follows that
\begin{equation}\label{FFrelations2RMM}
\phi(F_U(x))=F_X(x)\ \text{and}\  \psi(1-F_V(-y))=F_Y(y).
\end{equation}
For this example we assume that all three shocks are continuous random variables.

\begin{figure}[h]
\begin{center}
\includegraphics[width=9cm]{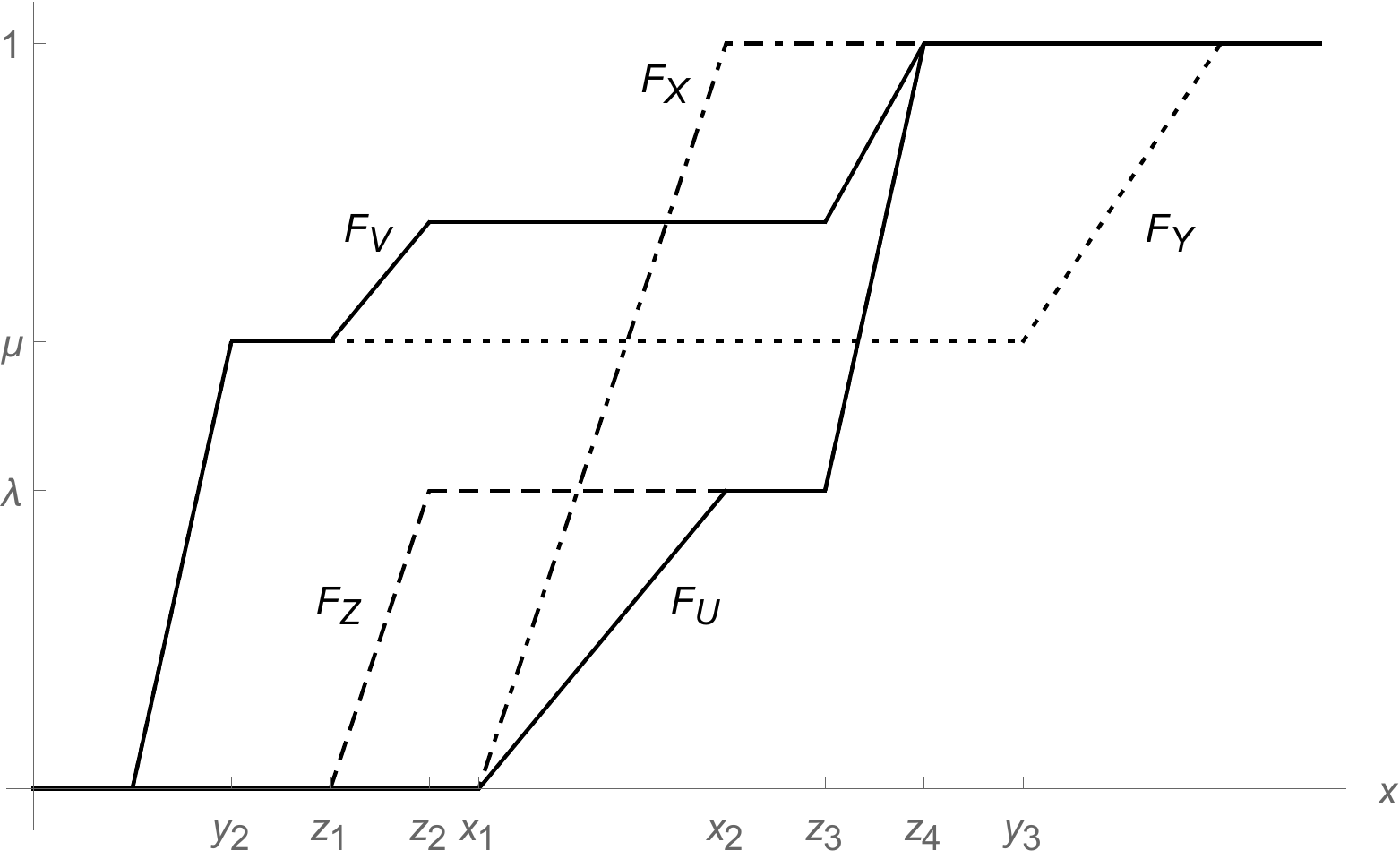}
\caption{Stochastic interpretation of copula $E_{\lambda,\mu}$.}\label{slika23}
\end{center}
\end{figure}

Consider now the reflected maxmin copula $E_{\lambda,\mu}$ of (\ref{d*_RMM}). Its generators are given by (\ref{f_lam g_mu}). Then we have
$$\phi_{\lambda}(t)=\left\{\begin{array}{ll}
        \frac{t}{\lambda}; & \text{if } 0\le t\le \lambda, \\
        1; & \text{if } t\ge \lambda,
        \end{array}\right. \ \text{and}\
\psi_{\mu}(t)=\left\{\begin{array}{ll}
        t; & \text{if } 0\le t\le 1-\mu, \\
        1-\mu; & \text{if }  1-\mu\le t<1,\\
        1; & \text{if } t=1.
        \end{array}\right. $$
The first relations of (\ref{FFrelationsRMM}) and (\ref{FFrelations2RMM}) imply that $\supp(Z)\subset (-\infty,z_2]\cup [z_3,\infty)$ and $\supp(X)\subset [x_1,x_2]$, where $z_2\le x_1\le x_2\le z_3$. We have that $F_Z(z_2)=P[Z\le z_2]=\lambda$ and $P[Z\ge z_3]=1-\lambda$. The second relations of (\ref{FFrelationsRMM}) and (\ref{FFrelations2RMM}) imply that $\supp(Y)\subset (-\infty,y_2]\cup [y_3,\infty)$ and $\supp(Z)\subset [z_1,z_4]$, where $y_2\le z_1\le z_4\le y_3$. We also see that $F_Y(y_2)=P[Y\le y_2]=\mu$ and $P[Y\ge y_3]=1-\mu$. Furthermore, we have that
$$F_U(x)=\left\{\begin{array}{ll}
        \lambda F_X(x); & \text{if } x\le z_3, \\
        F_Z(x); & \text{if } x\ge z_3,
        \end{array}\right.$$
and
$$F_V(y)=\left\{\begin{array}{ll}
        F_Y(y); & \text{if } y<z_1, \\
        \mu+(1-\mu) F_z(y); & \text{if } z_1\le y\le z_4,\\
        1; & \text{if } y\ge z_4.
        \end{array}\right.$$

Graphs of possible distribution functions of shocks $X$, $Y$ and $Z$ are illustrated in Figure \ref{slika23}. It clearly shows the interval $[x_1,x_2]$ where shock $X$ is acting, the intervals $[z_1,z_2]$ and $[z_3,z_4]$ where $Z$ is acting. These lie before and after $x_1$ and $x_2$, respectively. Finally, the shock $Y$ is acting partly before and partly after shock $Z$. Graphs of $F_U$ and $F_V$ are traced using a full line.

\section{Conclusion}

It is known since 1959 (Sklar's theorem) that given the marginal distributions every dependence in data can be described via a copula. However, it is not clear how to choose the family of copulas for our data in order to describe it at the best. We hope that this paper will help a practitioner to make a better choice based on the information about the non-exchangeability of the data.

Exchangeability as a probability concept (first next to independence) has been studied since 1930's; however, the point of non-exchangeability, also called asymmetry, of copulas was brought up only in 2006 and has been intensively studied ever since. One of the main contributions of this paper is an appropriate exertion of the \emph{maximal asymmetry function} on a particular family of copulas. This function was first introduced by Klement and Mesiar \cite{KM} and used on the family $\mathcal{C}$ of all copulas. It plays a major role in our paper, restricted to special families of copulas, due to the fact that it is computable for these families thus helping us determine the sharp bound of measures of asymmetry for each of the given families.

We are analyzing one of the most important classes of copulas for applications. Namely, we study the families of shock-based copulas, i.e.\ Marshall copulas, maxmin and reflected maxmin (RMM for short) copulas. They all belong to one of the larger families, either positive quadrant dependent or negative quadrant dependent copulas (respectively PQD or NQD for short), so we compute the function under consideration for these families as well. Using this gadget, we give the final results in measuring asymmetry of shock-based copulas. We compute the sharp bound of asymmetry measure $\mu_\infty$, the most important of the known asymmetry measures, for the family of Marshall's copulas and the family of maxmin copulas, which both equal to $\frac{4}{27}$~($\approx 0.148$). One should compare this bound to the one for the class of PQD copulas to which they belong, which is $3-2\sqrt{2}$~($\approx 0.172$), and to the general bound for all copulas that is $\frac13$. Furthermore, we give the sharp bound of the same asymmetry measure for RMM copulas which is $3-2\sqrt{2}$, compared to the same bound for NQD copulas, where they belong, which is $\sqrt{5}-2$~($\approx 0.236$).

The maximal asymmetry function is also useful in developing another major contribution of this paper. We give a detailed analysis of shocks in a given model at which the bound for asymmetry is attained. These interpretations for the three families studied are illustrated by examples that should be helpful to practitioners when choosing the model for their data.

%\newpage

\medskip

%\address
\noindent{\sc D. Kokol Bukov\v sek}: {Faculty of Economics, University of Ljubljana, and Institute of Mathematics, Physics and Mechanics, Ljubljana, Slovenia.}

%\email
\noindent{\it e-mail}: {\tt damjana.kokol.bukovsek@ef.uni-lj.si}%

\medskip
%\address
\noindent{\sc T. Ko\v sir}: {Faculty of Mathematics and Physics,
University of Ljubljana, and Institute of Mathematics, Physics and Mechanics,  Ljubljana, Slovenia.}

%\email
\noindent{\it e-mail}: {\tt tomaz.kosir@fmf.uni-lj.si}%

\medskip
%\address
\noindent{\sc B. Moj\v skerc}: {Faculty of Economics, University of Ljubljana, and Institute of Mathematics, Physics and Mechanics, Ljubljana, Slovenia.}

%\email
\noindent{\it e-mail}: {\tt blaz.mojskerc@ef.uni-lj.si}%

\medskip
%\address
\noindent{\sc M. Omladi\v c}: {Institute of Mathematics, Physics and Mechanics, Ljubljana, Slovenia.}

%\email
\noindent{\it e-mail}: {\tt matjaz@omladic.net}%

\vskip 7mm

\noindent{\large\bf Appendix: {Proofs of technical lemmas}}
\vskip 5mm

\noindent In the appendix proofs of Lemmas \ref{ocenaM}, \ref{ocenaMM} and \ref{ocenaIMM} are given. The underlying ideas of these proofs are similar, but the details of proofs significantly differ. Let us briefly describe the basic idea.

For given $x,y \in [0,1]$ and copula $C$ from a chosen family we estimate the difference $|C(x,y)-C(y,x)|$ with the use of one (Lemma \ref{ocenaM}) or two parameters (Lemmas \ref{ocenaMM} and \ref{ocenaIMM}) determined by the generating functions of copula $C$. Then we maximize the estimate over all possible values of these parameters.

\vspace{6mm}

\noindent {\bf \ref{ocenaM} Lemma.}\emph{ Let $C$ be a Marshall copula and $x, y \in [0, 1]$ with $y \ge x$. Then
$$|C(x, y) - C(y, x)| \le F(x,y)$$
where
$$F(x, y) = \left\{\begin{array}{ll}
        x(1-y); & \text{if } x \le y^2, \\
        \frac{x(y-x)}{y}; & \text{if } x > y^2. \\
        \end{array}\right.$$
}\\

\prf
If $C(x,y)=C(y,x)$ then the claim holds. So, suppose that $C(x,y)\neq C(y,x)$. In particular, this implies that $x\neq 0$, $y\neq 1,$ and $x\neq y$. We may assume without loss that
\begin{equation}\label{min<min}
\min\{xf(y),yg(x)\}<\min\{yf(x),xg(y)\}.
\end{equation}
If the latter inequality does not hold then we interchange the roles of $f$ and $g$.
\medskip

Now, we wish to show that $xf(y)\le yg(x)$. Assume to the contrary that $xf(y)> yg(x)$. Since $f^*$ and $g^*$ are nonincreasing and $y\ge x$ we have
$$\frac{g(y)}{y}\le\frac{g(x)}{x}<\frac{f(y)}{y}\le\frac{f(x)}{x}.$$
Then it follows that
$$\min\{yf(x),xg(y)\}=xg(y)<yg(x)=\min\{xf(y),yg(x)\},$$
which is in a contradiction with (\ref{min<min}). Therefore we have $xf(y)\le yg(x)$. We conclude that
\begin{equation}\label{Mar1}
|C(x,y)-C(y,x)|=%\min\{yf(x),xg(y)\}-\min\{xf(y),yg(x)\}=
\min\{yf(x),xg(y)\}-xf(y)\le\min\{yf(x),x\}-xf(y).
\end{equation}
The latter inequality follows since $g(y)\le 1$.
\medskip

We write $a=f(x)$. Since $\frac{f(x)}{x}\ge\frac{f(1)}{1}=1$ it follows that $a\ge x$, and thus $a\in [x,1]$. Also, we have $f(y)\ge\max\{y,a\}.$ We use these and (\ref{Mar1}) to obtain
\begin{equation*}\label{Mar2}
|C(x,y)-C(y,x)|\le \min\{ay,x\}-\max\{xy,xa\}.
\end{equation*}
It follows that we can take $F(x,y)=\max_{a\in[x,1]}\left\{\min\{ay,x\}-\max\{xy,xa\}\right\}.$
\medskip

Consider first the case $x\le y^2$. We separate further the cases $a\in [x,y]$ and $a\in [y,1]$. Then
\begin{eqnarray*}
\max_{a\in [x,y]}\{\min\{ay,x\}-\max\{xy,xa\}\}&=&\max_{a\in [x,y]}\{\min\{ay,x\}-xy\}\\
&=&\max_{a\in [x,y]}\{\min\{(a-x)y,x(1-y)\}\}\\
&=&\min\{(y-x)y,x(1-y)\}\\
&=&x(1-y).
\end{eqnarray*}
In the case $a\in\ [y,1]$ we have that $x(1-a)\le x(1-y)\le y^2-xy\le a(y-x)$. Thus, we obtain
\begin{eqnarray*}
\max_{a\in [y,1]}\{\min\{ay,x\}-\max\{xy,xa\}\}&=&\max_{a\in [y,1]}\{\min\{ay,x\}-ax\}\\
&=&\max_{a\in [y,1]}\{\min\{a(y-x),x(1-a)\}\}\\
&=&\max_{a\in [y,1]}\{x(1-a)\}\\
&=&x(1-y).
\end{eqnarray*}
Thus, we conclude that $F(x,y)=x(1-y)$ if $x\le y^2$.
\medskip

Next, consider the case $x\ge y^2$. Then
\begin{eqnarray*}
\max_{a\in [x,y]}\{\min\{ay,x\}-\max\{xy,xa\}\}&=&\max_{a\in [x,y]}\{\min\{(a-x)y,x(1-y)\}\}\\
&=&\min\{(y-x)y,x(1-y)\}\\
&=&(y-x)y
\end{eqnarray*}
and
\begin{eqnarray*}
\max_{a\in [y,1]}\{\min\{ay,x\}-\max\{xy,xa\}\}&=&\max_{a\in [y,1]}\{\min\{ay,x\}-ax\}\\
&=&\max_{a\in [y,1]}\{\min\{a(y-x),x(1-a)\}\}.
\end{eqnarray*}
Then, we observe that
$$\min\{a(y-x),x(1-a)\}=\left\{\begin{array}{ll} a(y-x);&a\le\frac{x}{y},\\ x(1-a);&a\ge\frac{x}{y}.\end{array}\right.$$
Maximum of this expression over all $a\in [y,1]$ is attained at $a=\frac{x}{y}$. Therefore, we have that
$$\max_{a\in[y,1]}\{a(y-x),x(1-a)\}=\frac{x(y-x)}{y}.$$ In this case, we conclude that $F(x,y)=\frac{x(y-x)}{y}$.
\qed

\vspace{6mm}

\noindent {\bf \ref{ocenaMM} Lemma.} \emph{Let $x, z \in [0, 1]$ and $x+z \le 1$. For $f,g\in\mF$ we have
$$|\min\{xz, f(x)g(z)\} - \min\{(1-x)(1-z), f(1-z)g(1-x)\}| \le G(x,z)$$
where
$$G(x, z) = \left\{\begin{array}{ll}
        \frac{x(1-x-z)}{1-z}; & \text{if } 1 - \sqrt{x} \le z \le \min\{x, 1-x\}, \\
        \frac{z(1-x-z)}{1-x}; & \text{if } 1 - \sqrt{z} \le x \le \min\{z, 1-z\}, \\
        xz; & \text{otherwise. }
        \end{array}\right.$$
}\\

\prf
Since $f(x)+x$ is a nondecreasing function on $[0,1]$ and $f(1)=0$, we have $f(x)+x\le f(1)+1=1,$ or $f(x)\le 1-x$ for all $x\in [0,1]$. Similarly, $g(x)\le 1-x$. Then $f(1-z)g(1-x)\le xz\le (1-x)(1-z)$ since $x+z\le 1$. This implies
\begin{equation}\label{min-min}
|\min\{ x z,f(x)g(z)\}-\min\{(1-x)(1-z),f(1-z)g(1-x)\}|=|\min\{ x z,f(x)g(z)\}-f(1-z)g(1-x)|.
\end{equation}
We write $a=f(1-z)$ and $b=g(1-x)$. Note that $0\le a\le z$ and $0\le b\le x$. Since $x\le 1-z$ and $z\le 1-x$ we have $f(x)+x\le f(1-z)+1-z$ and $g(z)+z\le g(1-x)+1-x$, respectively. Thus
\begin{equation}\label{min_min2}
f(x)\le 1+a-x-z \text{\ and\ }g(z)\le 1+b-x-z.
\end{equation}
In particular, we see that $f(x)\le a$ and $g(z)\le b$.
Since $\frac{f(x)}{x}$ and $\frac{g(z)}{z}$ are nonincreasing, we have $\frac{f(x)}{x}\ge \frac{f(1-z)}{1-z}$ and $\frac{g(z)}{z}\ge \frac{g(1-x)}{1-x}$. Hence
\begin{equation}\label{min_min3}
f(x)\ge \frac{ax}{1-z} \text{\ and\ }g(z)\ge \frac{bz}{1-x}.
\end{equation}
Now, we apply (\ref{min_min2}) and (\ref{min_min3}) to the right hand side of expression (\ref{min-min}) to obtain
\begin{equation}\label{min-min4}
\min\{ x z,f(x)g(z)\}-ab\le\min\{ x z-ab,(1+a-x-z)(1+a-x-z)-ab\}
\end{equation}
and
\begin{equation}\label{min-min5}
ab-\min\{ x z,f(x)g(z)\}\le ab-\min\left\{ x z,\frac{abxz}{(1-x)(1-z)}\right\}.
\end{equation}
We use the two bounds just proved to see that the expression in (\ref{min-min}) is bounded from above by
\begin{equation*}
\max\left\{\min\{ x z,(1+a-x-z)(1+b-x-z)\}-ab,ab-\min\left\{xz,\frac{abxz}{(1-x)(1-z)}\right\}\right\}.
\end{equation*}
Since $a\in [0,x]$ and $b\in [0,z]$, it follows that the expression in (\ref{min-min}) is bounded from above by
\begin{equation}\label{min-min6}
\max_{a\in [0,z],\,b\in [0,x]}\left\{\min\{ x z,(1+a-x-z)(1+b-x-z)\}-ab,ab-\min\left\{xz,\frac{abxz}{(1-x)(1-z)}\right\}\right\}.
\end{equation}
Note that $ab-xz\le 0$ and therefore we have
\begin{eqnarray*}\label{min-min7}
ab-\min\left\{xz,\frac{abxz}{(1-x)(1-z)}\right\}&=&\max\left\{ab-xz,\frac{ab(1-x-z)}{(1-x)(1-z)}\right\}\\
&=&\frac{ab(1-x-z)}{(1-x)(1-z)}\\
&\le&\frac{xz(1-x-z)}{(1-x)(1-z)}.
\end{eqnarray*}
We write $t=1-x-z$. Under our assumption $x+z\le 1$ we have $t\ge 0$. The first expression of the two that we maximize in (\ref{min-min6}) is then equal to
$$\min\{ x z-ab,(1+a-x-z)(1+b-x-z)-ab\}=\min\{ x z-ab,(t+a)(t+b)-ab\}.$$
We want to compute its maximal value over all $a\in [0,z]$ and $b\in [0,x]$. Since the expression $xz-ab$ is decreasing as $a$ and $b$ increase and the expression $(t+a)(t+b)-ab=t^2+at+bt$ is increasing as $a$ and $b$ increase their minimum is maximal at some point where they intersect. Thus, their minimum is attained on the curve $(t+a)(t+b)=xz$. In the case where $xz\le t^2$, we have $xz\le (t+a)(t+b)$ and
\begin{equation*}
\max_{a\in [0,z],\,b\in [0,x]}\min\{ x z-ab,(t+a)(t+b)-ab\}=xz.
\end{equation*}
Suppose next that $xz\ge t^2$. In this case, there are solutions in $a$ and $b$ of the equation $xz=(t+a)(t+b)$ that lie in the rectangle $[0,z]\times[0,x]$. Then
\begin{equation*}
\max_{a\in [0,z],\,b\in [0,x]}\min\{ x z-ab,(t+a)(t+b)-ab\}=\max_{a\in [0,z],\,b\in [0,x],\,xz=(t+a)(t+b)}\{(a+b)t+t^2\}.
\end{equation*}
To find this maximum, we have to analyze several cases with respect to $x$ and $z$. The maximum is attained on the boundary of the rectangle, i.e., at points $(0,\frac{xz-t^2}{t})$ or  $(\frac{xz-t^2}{t},0)$ or  $(z,\frac{z(x-t)-t^2}{t+z})$ or  $(\frac{x(z-t)-t^2}{t+x},x)$, whichever is greater and possible. The case $a=z$ is possible only if $b=\frac{z(x-t)-t^2}{t+z}\ge 0$, or equivalently, only if $z\ge (1-x)^2$. In this case we have $(a+b)t+t^2=\frac{zt}{t+z}=\frac{z(1-x-z)}{1-x}$. The case $b=x$ is possible only if $a=\frac{x(z-t)-t^2}{t+x}\ge 0$, or equivalently, only if $x\ge (1-z)^2$. In this case, we have $(a+b)t+t^2=\frac{xt}{t+x}=\frac{x(1-x-z)}{1-z}$. In the cases $a=0$ or $b=0$, we have $(a+b)t+t^2=xz.$ Since $x+z\le 1$, it follows that $\frac{z(1-x-z)}{1-x}\ge \frac{x(1-x-z)}{1-z}$ if and only if $z\ge x$.
Then, the above analysis of cases implies that (\ref{min-min}) is bounded from above by
\begin{equation*}
 \left\{\begin{array}{ll}
        \frac{x(1-x-z)}{1-z}; & \text{if } x\ge z \text{ and }x\ge(1-z)^2, \\
        \frac{z(1-x-z)}{1-x}; & \text{if } z\ge x \text{ and }z\ge(1-x)^2,  \\
        xz; & \text{if } x\le (1-z)^2 \text{ or }z\le(1-x)^2.
        \end{array}\right.
\end{equation*}
Now, we see that $\frac{xz(1-x-z)}{(1-x)(1-z)}\le G(x,z)$ for any $x,z\ge 0$ with $x+z\le 1$, so the maximum in expression (\ref{min-min6}) is equal to $G(x,z)$, which completes the proof.
\qed

\vspace{6mm}

\noindent {\bf \ref{ocenaIMM} Lemma.}  \emph{
Let $x, y \in [0, 1]$, $y \ge x$ and suppose $f(x)g(y) \le f(y)g(x)$. Then
$$|\max\{0, xy - f(x)g(y)\} - \max\{0, xy - f(y)g(x)\}| \le H(x,y)$$
where
$$H(x, y) = \left\{\begin{array}{ll}
        \frac{x(y-x)}{y}; & \text{if } y(1-y) \le x \le 1-y, \\
        \frac{(y-x)(1-y)}{y}; & \text{if } x \ge 1-y, \\
        xy; & \text{if } x \le y(1-y).
        \end{array}\right.$$
}\\

\prf
We may and do assume that $x\neq 0$, $y\neq 1$ and $x\neq y$. Since $y \ge x$ and  $f(y)g(x)\ge f(x)g(y)$ we have
\begin{equation}\label{imaxmin}
|\max\{0, xy - f(x)g(y)\} - \max\{0, xy - f(y)g(x)\}| = \max\{0, xy - f(x)g(y)\} - \max\{0, xy - f(y)g(x)\}.
\end{equation}
Now either both maxima in (\ref{imaxmin}) are zero or $xy-f(x)g(y)$ is positive. Assume that the latter holds. We write $a=f(x)$ and $b=g(y)$. Then $xy\ge ab$. Using these we see that (\ref{imaxmin}) is equal to
\begin{eqnarray*}\label{imaxmin2}
&=& xy-f(x)g(y)-\max\{0,xy-f(y)g(x)\}\\
&=& \min\{xy-f(x)g(y),f(y)g(x)-f(x)g(y)\}\\
&=& \min\{xy-ab,f(y)g(x)-ab\}.
\end{eqnarray*}
Recall that $0\le f(x)\le 1-x$ for each function $f\in\mF$. So, we have $a\in [0,1-x]$, $b\in[0,1-y]$ and $ab\le xy$.

Since function $g(x)+x$ is nondecreasing and $x\le y$, it follows that $g(x)+x\le g(y)+y$ and so $g(x)\le b+y-x$.

Since function $f^*(x)=\frac{f(x)}{x}$ is nonincreasing, we have $\frac{f(y)}{y}\le\frac{f(x)}{x}$, and thus $f(y)\le \frac{ay}{x}$. Combining these relations, we get $f(y)\le\min\{\frac{ay}{x},1-y\}$.

Using all of the above, we see that the expression in (\ref{imaxmin}) is bounded above by
\begin{eqnarray}
&\le& \min\left\{xy-ab,\frac{ay}{x}(b+y-x)-ab,(1-y)(b+y-x)-ab\right\}\\
&\le& \max_{a\in[0,1-x],b\in[0,1-y],ab\le xy}\min\left\{xy-ab,\frac{(y-x)}{x}a(b+y),(1-y)(b+y-x)-ab\right\}. \label{imaxmin3}
\end{eqnarray}
We denote the latter expression by $H(x,y)$.

Suppose first that $x+y\ge 1$. Then $a\le 1-x$ and $b\le 1-y$ imply $ab\le 1-x-y+xy\le xy$ and $(1-y)(b+y-x)\le (1-y)(1-x)\le xy$. So the expression in (\ref{imaxmin3}) is equal to
\begin{equation}\label{imaxmin4}
\max_{a\in[0,1-x],b\in[0,1-y]}\min\left\{\frac{(y-x)}{x}a(b+y),(1-y)(b+y-x)-ab\right\}.
\end{equation}
The first function as a function of $a$ and $b$ inside braces in (\ref{imaxmin4}) is increasing in $a$, while the second function is decreasing in $a$. So the maximal value over $a$ is attained at the intersection of graphs of the two functions. There we have $\frac{y-x}{x}a(b+y)=(1-y)(b+y-x)-ab$. Since $b+y-x>0$ it follows that $a=\frac{x(1-y)}{y}$, which lies in the interval $[0,1-x]$ as required. Then, we have
$$H(x,y)=\max_{b\in [0,1-y]}\left\{\frac{(y-x)(1-y)(b+y)}{y}\right\}=\frac{(y-x)(1-y)}{y}.$$

Suppose now that $x+y\le 1$ and consider first the case $x\le y(1-y)$. Since $b\ge 0$ we have that
\begin{equation}\label{imaxmin5}
(1-y)(b+y-x)\ge (1-y)(y-x)=y(1-y)-x+xy\ge xy.
\end{equation}
Therefore, we can omit the last term in expression (\ref{imaxmin3}) for $H(x,y)$ and we obtain
\begin{equation}\label{imaxmin6}
H(x,y)=\max_{a\in[0,1-x],b\in[0,1-y],ab\le xy}\min\left\{xy-ab,\frac{(y-x)}{x}a(b+y)\right\}.
\end{equation}
The first function as a function of $a$ and $b$ inside braces in (\ref{imaxmin6}) is decreasing in $a$, while the second function is increasing in $a$. Thus, the maximal value over $a$ is attained at the intersection of graphs of the two functions. There we have $x^2y={ay}(b+y-x)$, and so $a=\frac{x^2}{b+y-x}$. We use (\ref{imaxmin5}) and assumption $x\le y$ to show that
\begin{equation*}
a=\frac{x^2}{b+y-x} \le \frac{(1-y)x}{y}\le 1-x.
\end{equation*}
We also have
\begin{equation}\label{imaxmin7}
xy-ab=xy-\frac{x^2b}{b+y-x}=\frac{x(b+y)(y-x)}{b+y-x}\ge 0.
\end{equation}
Then it follows using (\ref{imaxmin7}) that
\begin{equation*}
H(x,y)=\max_{b\in [0,1-y]}\left\{\frac{x(b+y)(y-x)}{b+y-x}\right\}.
\end{equation*}
This is a rational function in $b$ with zero at $b=-y$, pole at $b=x-y$ and horizontal asymptote at $x(y-x)$. Therefore, it is decreasing on the interval $[0,1-y]$, so the maximum is attained at $b=0$. Thus we have $H(x,y)=xy$ in this case.

Suppose finally that $y(1-y)\le x \le 1-y$. We consider for a moment two separate cases. In the case that $b+y\le\frac{x}{1-y}$ the last term in expression (\ref{imaxmin3}) is smaller than the first term. Then we have
\begin{equation}\label{imaxmin8}
H(x,y)=\max_{a\in[0,1-x],b\in\left[0,\frac{x-y+y^2}{1-y}\right]}\min\left\{\frac{(y-x)}{x}a(b+y),(1-y)(b+y-x)-ab\right\}.
\end{equation}
Similar arguments as were used in the previous cases imply that the maximum in this case is attained at $a=\frac{x(1-y)}{y}$.
In the other case when $b+y\ge\frac{x}{1-y}$, we see that the first term in expression (\ref{imaxmin3}) is smaller than the last term. Then we have
\begin{equation}\label{imaxmin9}
H(x,y)=\max_{a\in[0,1-x],b\in\left[\frac{x-y+y^2}{1-y},1-y\right]}\min\left\{xy-ab,\frac{(y-x)}{x}a(b+y)\right\}.
\end{equation}
Similarly as before, we show that the maximum in this case is attained at $a=\frac{x^2}{b+y-x}$. Combining the two cases in (\ref{imaxmin8}) and (\ref{imaxmin9}) we obtain
$$
H(x,y)=\max\left\{\max_{b\in\left[0,\frac{x-y+y^2}{1-y}\right]}\left\{\frac{(y-x)(1-y)(b+y)}{y}\right\}, \max_{b\in\left[\frac{x-y+y^2}{1-y},1-y\right]}\left\{\frac{x(y-x)(b+y)}{b+y-x}\right\}\right\}.
$$
The first maximum inside the braces is attained when $b$ is as large as possible, while the second maximum is attained when $b$ is as small as possible. So we have $b=\frac{x-y+y^2}{1-y}$ and $H(x,y)=\frac{x(y-x)}{y}$. The proof is complete.
\qed

\end{document}